\documentclass[]{interact}

\usepackage{epstopdf}
\usepackage[caption=false]{subfig}

\usepackage[numbers,sort&compress]{natbib}
\bibpunct[, ]{[}{]}{,}{n}{,}{,}

\theoremstyle{plain}
\newtheorem{theorem}{Theorem}[section]
\newtheorem{lemma}[theorem]{Lemma}
\newtheorem{corollary}[theorem]{Corollary}

\theoremstyle{definition}
\newtheorem{definition}[theorem]{Definition}

\theoremstyle{remark}

\usepackage{algorithm}\usepackage{algpseudocode}
\newcommand{\argmin}{\mathop{\arg\!\min}}

\usepackage{nicefrac}
\usepackage{xcolor}

\allowdisplaybreaks

\def \R {\mathbb R}
\def\eqdef{\overset{\text{def}}{=}}

\newcommand{\EndProof}{\begin{flushright}$\square$\end{flushright}}

\newcommand{\cN}{{\cal N}}

\newcommand{\EE}{\mathbf{E}}

\def\R{\mathbb{R}}
\newcommand{\E}{{\mathbb E}}

\def\X{\mathcal X}
\def\Y{\mathcal Y}

\def\R{\mathbb R}
\def\E{\mathbb E}
\def\EE{\mathbb E}

\def\la{\langle}
\def\ra{\rangle}

\def\one{{\mathbf 1}}

\newcommand{\ab}[1]{{\color{black}#1}}
\newcommand{\abc}[1]{{\color{black}#1}}
\newcommand{\abcd}[1]{{\color{black}#1}}

\begin{document}


\title{One-Point Feedback for Composite Optimization with Applications to Distributed and Federated Learning}

\author{
\name{Aleksandr Beznosikov\textsuperscript{a,b,c}\thanks{CONTACT Aleksandr Beznosikov. Email: anbeznosikov@gmail.com}, Ivan Stepanov\textsuperscript{b,a}, Artyom Voronov\textsuperscript{a} and Alexander Gasnikov\textsuperscript{c,a,b}}
\affil{\textsuperscript{a}Moscow Institute of Physics and Technology, Moscow, Russia; \textsuperscript{b}Ivannikov Institute for System Programming, Moscow, Russia; \textsuperscript{c}Innopolis University, Innopolis, Russia}
}

\maketitle

\begin{abstract}
This work is devoted to solving the composite optimization problem with the mixture oracle: for the smooth part of the problem, we have access to the gradient, and for the non-smooth part, only the one-point zero-order oracle is available. For such a setup, we present a new method based on the sliding algorithm. Our method allows to separate the oracle complexities and to compute the gradient for one of the functions as rarely as possible. The paper also presents the applicability of our new method to the problems of distributed optimization and federated learning. Experimental results confirm the theory.
\end{abstract}


\section{Introduction}

\textbf{Composite optimization.} In this paper, we focus on the composite optimization problem \cite{nesterov2018lectures,lan}:
\begin{equation}
    \label{problem_orig}
    \min_{x \in \mathcal{X}} \Psi_0 (x) \eqdef f(x) + g(x),
\end{equation}
\ab{where we minimize the two-part function $\Psi_0$ on some set $\X$.}
This problem occurs in a fairly large number of applications. In particular, we can recall the problems of minimizing the objective function $f(x)$ with regularization $g(x)$, which can often be found in machine learning \cite{bach2021learning}. 
Newer and more interesting applications of the composite problem arise in distributed optimization. In more details, the goal of distributed optimization is to minimize the global objective function $f(x) = \sum_{m=1}^M f_m (x)$, where functions $f_1, \ldots, f_M$ are distributed over $M$ devices/workers, and each device $m$ has access only to its local function $f_m$.  Therefore, in order to solve this problem, one needs to establish a communication process between the devices. There are two methods: centralized and decentralized. In the centralized case, all devices can communicate only with the central server -- transfer information about the $f_m$ function to it and receive responses. In the decentralized setting, there is no central server; all devices are connected into a network, which can be represented as an undirected graph, \abc{where vertices represent devices, and edges represent the presence of the connection between a pair of devices.} Communication in the decentralized network is typically done with the gossip protocol \cite{kempe2003gossip,boyd2006randomized,nedic2009distributed}, which uses the so-called gossip matrix \ab{$W \in \R^{M \times M}$}. This matrix is built on the basis of the properties of the communication graph. 

It turns out that
the distributed optimization problem can be written as a composite one \cite{9173772,gorbunov2019optimal,beznosikov2019derivative,dvinskikh2021decentralized,hanzely2021personalized,hanzely2020lower}:
\begin{equation}
    \label{centr}
    \min_{(x_1, \ldots, x_M) \in \mathcal{X}^M} \underbrace{\frac{1}{M}\sum\limits_{m=1}^M f_m(x_m)}_{f(x_1, \ldots, x_M)} + \underbrace{\frac{\lambda}{2 M}\sum\limits_{m=1}^M \|x_m - \bar x\|^2_2}_{g(x_1, \ldots, x_M)}
\end{equation}
for the centralized case and
\begin{equation}
    \label{decentr}
    \min_{(x_1, \ldots, x_M) \in \mathcal{X}^M} \underbrace{\frac{1}{M}\sum\limits_{m=1}^M f_m(x_m)}_{f(x_1, \ldots, x_M)} + \underbrace{\frac{\lambda}{2 M} \|\sqrt{W} X\|^2_F}_{g(x_1, \ldots, x_M)}
\end{equation}
for decentralized one. Here we introduce the matrix $X = [x_1, \ldots, x_M]^T$, the vector $\bar x = \frac{1}{M} \sum_{m=1}^M x_m$ and the parameter of regularization $\lambda> 0$. The essence of the expressions \eqref{centr} and \eqref{decentr} are very simple. On each device, we have local variables $x_m$, and we penalize their deviations at the expense of the regularizer. In the centralized case, we penalize the deviation from the average across the entire network, and in the decentralized case, the difference between the connected devices (this is what the $W$ matrix is responsible for). In fact, in the decentralized case, we can also write the penalized problem in form \eqref{centr}. \abc{\abcd{However,} while in the case of a centralized architecture $\bar x$ is easy to compute on the server, in a decentralized network this is problematic (in particular, one of the devices has to be used as a server).} Another important question is how to choose \abc{the parameter $\lambda$}. \abcd{To get a solution close to the optimal one} of the distributed problem, one needs to take $\lambda$ large enough \cite{dvinskikh2021decentralized}. But more recently, the problems \eqref{centr} and \eqref{decentr} were considered from the point of view of personalized federated learning, in which case it makes sense to take small $\lambda$ as well \cite{hanzely2020lower,hanzely2021personalized,sadiev2021decentralized}.

\textbf{Gradient-free methods.} Now let us go back to the original problem \eqref{problem_orig}. As noted above, the function $g$ often plays the role of a regularizer, usually it is a simple function for which the gradient $\nabla g$ can be computed. At the same time, the objective function $f$ can be quite complex. In this paper, we focus on the case when for the function $f$ only zero-order oracle 
(i.e., only the values of the function $f$, but not its gradient) is available. In the literature, this concept is sometimes referred to as a black box.  It arises when the calculation of gradient is expensive (in adversarial training \cite{Chen_2017}, optimization \cite{Nesterov},  structured-prediction learning \cite{Taskar04learningstructured}) or  impossible (in reinforcement learning \cite{fazel2018global,pmlr-v80-choromanski18a,salimans2017evolution}, bandit problem \cite{bubeck2012regret,shamir2017optimal}, black-box ensemble learning \cite{lian2015asynchronous}). To make the problem statement even more practical we assume that we have access to inexact values of the function $f(x, \xi)$ with some random noise $\xi$. With the help of this oracle, it is possible to make some approximation of the gradient in terms of finite differences. Next we highlight two main approaches for such gradient estimation. The first approach is called a two-point feedback:
\begin{equation}
\label{TPF}
\tilde f'_r(x,\xi) \eqdef \frac{n}{2r} (f(x+  r e, \xi)  - f(x - re,\xi)) e,
\end{equation}
where $e \in \R^n$ is uniformly distributed on the unit Euclidean sphere, \ab{$r > 0$ is the smoothing parameter responsible for the margin length of the approximation.} For the two-point feedback there \abc{are} many papers with theoretical analysis \cite{duchi2015optimal,Nesterov,gasnikov2017stochastic,shamir2017optimal,dvurechensky2021accelerated,gorbunov2022accelerated}.
An important thing about this approach is the assumption that we are able to obtain the values of the function at points $x+  r e$ and $x - r e$ with the same realization of the noise $\xi$. However, from a practical point of view, this is a very strong and idealistic assumption. Therefore, it is proposed to consider the concept of one-point feedback \cite{gasnikov2017stochastic,akhavan2020exploiting,novitskii2021improved}:
\begin{equation}
\label{grad_f}
\tilde f'_r(x,\xi^{\pm}) \eqdef \frac{n}{2r} (f(x+  r e, \xi^+)  - f(x - r e,\xi^-)) e .
\end{equation}
In general $\xi^+ \neq \xi^-$.  In this paper we work with the one-point concept.

The function $f$ is "bad", while the function $g$ is "good". The question arises how to minimize $\Psi_0$ from \eqref{problem_orig}. The easiest option is to add the gradient of $g$ and the "gradient" of $f$ (from \eqref{grad_f}) and make step along it. In this approach, there are no problems when $g$ is just a Tikhonov regularizer, but if we look at the problems \eqref{centr} and \eqref{decentr}, to compute the gradient of $g$ we need to make communication, while to calculate the "gradient" of $f$ we do not need it. But communications are the bottleneck of distributed algorithms, they require significantly more time than local computations. Therefore, one wants to reduce the number of communications and to calculate the gradient $g$ as rarely as possible. 

This brings us to the goal of this paper: to come up with an algorithm that solves the composite optimization problem for one part of which we have a one-point zero-order oracle, and for the other -- a gradient. At the same time, \abcd{we want to make calls to the gradient as rarely as possible.}

\subsection{Our contribution} 

We present a new method based on the sliding technique for the convex problem \eqref{problem_orig} with the mixture oracle: first-order for the smooth part $g$ and zero-order for the non-smooth part $f$.  Our method solves the problems mentioned above in the introduction. It reduces the number of calls \abcd{to} the gradient $\nabla g$ of the smooth part of the composite problem, while using one-point feedback for the non-smooth part $f$.

Note that all the results were obtained in the general (non-Euclidean) proximal setting \ab{to take into account the geometry of the problem. This was achieved by using the Bregman divergence instead of the classical Euclidean distance.} It allows sometimes to reduce the $f$-oracle calls complexity $\sim n$-times in comparison with algorithms that use the Euclidean setup, where $n$ is a dimension of the problem -- see Table~\ref{tab:summary1}.

\renewcommand{\arraystretch}{1.7}
\renewcommand{\tabcolsep}{6pt} 
\begin{table}[h!]
\begin{center}
\begin{tabular}{c|c|c}
\hline
 & $\|\cdot\|=\|\cdot\|_2$ & $\|\cdot\|=\|\cdot\|_1$  \\ \hline
$\nabla g$ & $O\left(\sqrt{\frac{L D^2}{\varepsilon}} \right)$ & $O\left(\sqrt{\frac{L D^2}{\varepsilon}} \right)$  \\ \hline
$f$ & $\tilde O\left(\sqrt{\frac{L D^2}{\varepsilon}} + \frac{n  G^2 D^2}{\varepsilon^2} + \frac{n^2  G^2\sigma^2 D^2}{\varepsilon^4} \right)$ & $O\left(\sqrt{\frac{L D^2}{\varepsilon}} + \frac{\ln n G^2 D^2}{\varepsilon^2} + \frac{n \ln n G^2\sigma^2 D^2}{\varepsilon^4} \right)$ \\ \hline
\end{tabular}
\caption{\small Summary of complexity results on the number of $\nabla g$ and $f$ calls for finding an $\varepsilon$-solution of \eqref{problem_orig} in the different geometric setups: Euclidean and non-Euclidean. Convergence is measured by the function distance to the solution.   {\em Notation:} $L$ = constant of $L$-smoothness of $g$ in the corresponding norm, \ab{$G$ = upper bound of $f'$ in the $\ell_2$-norm ($\| f'(x) \|_2 \le G$ for all $f'(x) \in \partial f(x)$ \abc{and $x \in \mathcal{X}$})}, $n$ = dimension of vector $x$ ($\mathcal{X} \subset \R^n$), $\sigma^2$ = variance of unbiased noise $\xi$, $D$ = diameter of $\mathcal{X}$(in the corresponding geometry).}
\label{tab:summary1}
\end{center}
\end{table}

We also present the applicability and relevance of our new method for distributed and federated learning problems in both centralized \eqref{centr} and decentralized \eqref{decentr} setups -- see Table~\ref{tab:summary2}. It turns out that this method can be useful in terms of reducing the number of communications.

\begin{table}[h!]
\vspace{-0.3cm}
\begin{center}
\begin{tabular}{c|c|c}
\hline
& \textbf{Centralized} & \textbf{Decentralized}  \\ \hline
{\tt comm}& $O\left(\sqrt{\frac{\lambda D^2}{\varepsilon}} \right)$ & $O\left(\sqrt{\frac{\lambda \lambda_{\max}(W) D^2}{\varepsilon}} \right)$  \\ \hline
{\tt local}&\!$O\!\left(\!\sqrt{\frac{\lambda D^2}{\varepsilon}}\!+\!\frac{M n  G^2 D^2}{\varepsilon^2}\!+\!\frac{M n^2  G^2\sigma^2 D^2}{\varepsilon^4} \!\right)$\! & $O\!\left(\!\sqrt{\frac{\lambda \lambda_{\max}(W) D^2}{\varepsilon}}\!+\!\frac{M n G^2 D^2}{\varepsilon^2}\!+\!\frac{M n^2 G^2\sigma^2 D^2}{\varepsilon^4}\!\right)$\! \\ \hline
\end{tabular}
\caption{\small Summary of complexity results on communications ({\tt comm}) and local computations ({\tt local}) for finding an $\varepsilon$-solution of centralized \eqref{centr} and decentralized \eqref{decentr} distributed problems. Convergence is measured by the function distance to the solution. {\em Notation:} $\lambda_{\max}(W)$ = maximum eigenvalue of $W$, \ab{$G$ = upper bound of $f'_m$ in the $\ell_2$-norm ($\| f'_m(x) \|_2 \le G$ for all $f'(x) \in \partial f(x)$ \abc{, $x \in \mathcal{X}$} and $m$)}, $n$ = dimension of vector $x$ ($\mathcal{X} \subset \R^n$), $\sigma^2$ = variance of unbiased noise $\xi$, $D$ = diameter of $\mathcal{X}$ (in the Euclidean norm).}
\label{tab:summary2}
\end{center}
\end{table}

\subsection{Comparison with known results}

Let us note some works related to our paper.

\textbf{Sliding.} 
The naive approach to \eqref{problem_orig} looks at it as a whole problem and does not take into account its composite structure. This can significantly worsen the oracle complexity \ab{(number of the oracle calls)} for one of the functions. The sliding technique allows to avoid these losses and to separate oracle complexities. In particular, if we can solve a separate problem $\min f(x)$ by $T_f$
oracle calls (these can be calls of gradient or any other oracle, for example, zero-order), and a problem $\min g(x)$ by $T_g$
oracle calls, then the sliding technique gives that we can solve the composite problem by $\mathcal{O}(T_f)$ oracle calls corresponding to $f$ and $\mathcal{O}(T_g)$ oracle calls corresponding to $g$. As mentioned above, if we use the naive approach we have the same complexity $\mathcal{O}(\max(T_g; T_f))$ for both $f$ and $g$.

There are various types of sliding in the literature, depending on what assumptions are made for \eqref{problem_orig}.
\begin{itemize}

\item  The sliding was justified for smooth $f$ and  $g$ with gradient oracles for convex optimization problems in \cite{juditsky2011solving,alkousa2020accelerated,vladislav2021accelerated,lan2021mirror,lin2022inexact}.

\item In \cite{dvinskikh2020accelerated,ivanova2020oracle}, they considered case with smooth $f$ and $g$ with zero-order oracle for $f$ and gradient oracle for $g$.

\item  In \cite{lan2016gradient,lan}, the sliding technique was used for non-smooth $f$ and smooth $g$ with stochastic subgradient oracle for $f$ and gradient oracle for $g$.

\item \abc{The authors of \cite{lan2016gradient}, \cite{beznosikov2019derivative} adapted} the sliding approach for non-smooth $f$ and smooth $g$ with two-point zero-order stochastic oracle \eqref{TPF} for $f$ and gradient oracle for $g$.
\end{itemize}
The development of the sliding technique is a quite popular issue in the literature, \abcd{but on the other hand that there are still many open problems especially for the mixture oracle.} In this paper, we concentrate on the generalization of \cite{beznosikov2019derivative} for non-smooth $f$ and smooth $g$ with \textbf{one}-point zero-order stochastic oracle for $f$ (one-point feedback rather than two-point of \cite{beznosikov2019derivative}) and a gradient oracle for $g$ for convex optimization problems. For strongly convex problems our results can also be generalized by using the standard restart technique, see e.g. \cite{beznosikov2019derivative}.



\textbf{Gradient-free methods.} Let us highlight the main works devoted to the zero-order methods: for two-point feedback \cite{Shamir15,Nesterov,duchi2015optimal,gasnikov2017stochastic,shamir2017optimal,dvurechensky2021accelerated,gorbunov2022accelerated}, for one-point feedback \cite{bach2016highly,gasnikov2017stochastic,akhavan2020exploiting,novitskii2021improved}.
For two-point stochastic/deterministic feedback optimal methods for smooth/non-smooth, convex/strongly convex problems were developed in the cited papers. For the one-point feedback setting, there \abc{is} still a gap between lower bounds and the complexities of the best known methods. In this paper, we generalize the best-known composite-free ($g = 0$) results concerning non-smooth $f$ with stochastic one-point feedback from \cite{gasnikov2017stochastic} for problems with smooth regularizer $g\neq 0$.

\textbf{Distributed setup.} For (strongly) convex optimization problems optimal (stochastic) gradient decentralized methods were developed -- see surveys \cite{dvinskikh2021decentralized,gorbunov2020recent} and references therein. For stochastic two-point feedback (with non-smooth target function $f$) optimal decentralized methods were developed in \cite{beznosikov2019derivative}. To the best of our knowledge this is the only optimal result in this field. For one-point stochastic feedback we know only one result \cite{akhavan2021distributed}, they assume that the target function $f$ is highly-smooth and strongly convex. \abcd{However,} the method is very expensive in terms of decentralized communications. The reason for this issue is that, in \cite{akhavan2021distributed}, \ab{the} authors fight only for the oracle calls criteria and do not use the sliding technique that allows to split communication complexity from the oracle one. \abc{In our paper, using the sliding technique, we split these complexities in the problem \eqref{decentr} and obtain much better guarantees on the number of communications.}


\section{Preliminaries}

First, we define some notation. We denote the inner product of the vectors $x,y\in\R^n$ as $\langle x,y \rangle \eqdef \sum_{i=1}^nx_i y_i$, where $x_i$ corresponds to the $i$-th component of $x$ in the standard basis in $\R^n$. We also define $\ell_p$-norms as $\|x\|_p \eqdef \left(\sum_{i=1}^n|x_i|^p\right)^{\frac{1}{p}}$ for $p\in(1,\infty)$ and for $p = \infty$ we denote $\|x\|_\infty \eqdef \max_{1\le i\le n}|x_i|$. \ab{We use $\| \cdot \| = \| \cdot \|_p$ for brevity.} The dual norm $\|\cdot\|_*$ for the norm $\|\cdot\|$ is defined as follows: $\|y\|_* \eqdef \max\left\{\langle x, y \rangle \mid \|x\| \le 1\right\}$. 
We use $q$ as the degree of the dual norm $\|\cdot\|_* = \| \cdot \|_q$. It is known that $p$ and $q$ are related by the following proportion: $\frac{1}{p} + \frac{1}{q} = 1$. 
\abc{We also introduce the norms of the real matrix $X \in \R^{m \times n}$: the spectral norm as $\|X\|_{2}= \lambda_{\max} (X^T X)$ and the Frobenius norm as $\|X\|_{F}={\sqrt {\sum _{i=1}^{m}\sum _{j=1}^{n}|x_{ij}|^{2}}}$. To define the Kronecker product of $X \in \R^{m \times m}$ and $B \in \R^{n \times n}$ we use $A \otimes B$. We denote the Minkowski sum of the sets $\X$ and $\Y$ as $\X + \Y$.}
The operator $\E[\cdot]$ denotes the full mathematical expectation and the operator $\E_\xi[\cdot]$ expresses the conditional mathematical expectation with respect to all randomness coming from the random variable $\xi$.

Now let us introduce a few definitions.
\begin{definition}[$L$-smoothness]
\label{def:smooth}
    The function $g$ is called $L$-smooth w.r.t. norm $\|\cdot\|$ on $\X \subset \R^n$ with $L > 0$ when it is differentiable and its gradient is $L$-Lipschitz continuous on $\X$, i.e.\ 
    \begin{equation*}
        \|\nabla g(x) - \nabla g(y)\|_* \le L\|x - y\|~~~\text{for all}~~~x,y\in \X.
    \end{equation*}
\end{definition}
One can show that $L$-smoothness implies \cite{nesterov2018lectures}
\begin{equation}
    \label{g-L-smooth} 
    g(x) \leq g(y) + \langle \nabla g(y), x - y\rangle + \dfrac{L}{2}\|x-y\|^2~~~\text{for all}~~~x,y\in \X.
\end{equation}
\begin{definition}[$G$-Lipschitzness]
    The function $f$ is called $G$-Lipschitz w.r.t. norm $\|\cdot\|$ on $\X \subset \R^n$ with $G > 0$ when it holds that
    \begin{equation*}
        |f(x) - f (y)| \le G\|x - y\|~~~\text{for all}~~~x,y\in \X.
    \end{equation*}
\end{definition}
\ab{
\begin{definition}[Convexity]
\label{def:conv}
    Continuously differentiable function $g$ is called $\mu$-strongly convex w.r.t. norm $\|\cdot\|$ on $\X \subset \R^n$ \ab{with $\mu > 0$} when it holds that
    \begin{equation*}
    g(y) \geq g(x) + \langle \nabla g(x), y - x \rangle + \frac{\mu}{2}\| x - y\|^2~~~\text{for all}~~~x,y\in \X.
    \end{equation*}
For a non-differentiable function $f$ the $\mu$-strong convexity w.r.t. norm $\|\cdot\|$ on $\X$ is introduced as follows:
\begin{equation*}
    f(\alpha x + (1 - \alpha) y) \leq \alpha f(x) + (1 - \alpha) f(y) -  \alpha (1 - \alpha)\frac{\mu}{2}\| x - y\|^2~~~\text{for all}~~~x,y\in \X, ~\alpha\in [0;1].
\end{equation*}
If $\mu = 0$, then $g$ and $f$ are convex.
\end{definition}
}
\begin{definition}[Bregman divergence]
    Suppose some continuously differentiable function $\nu(x)$ is $1$-strongly convex w.r.t. norm $\|\cdot\|$ on $\X$. Then for any two points $x,y\in \X$ we define the Bregman divergence $V(x,y)$ associated with $\nu(x)$ as follows:
    \begin{equation*}
        V(x,y) \eqdef  \nu(y) - \nu(x) -  \langle \nabla \nu(x), y-x \rangle.
    \end{equation*}
We denote the Bregman diameter of the set $\X$ w.r.t.\ $V(x,y)$ as $D_{\X,V} \eqdef \max\{\sqrt{2V(x,y)}\mid x,y \in \X\}$.
\end{definition}

\section{Main part}

Recall that we consider the composite optimization problem \eqref{problem_orig}. To take into account the "geometry" of the problem, we work in a certain (not necessarily Euclidean) norm $\| \cdot \|$ (with dual norm $\| \cdot \|_*$), and also measure the distance using the Bregman divergence $V$. Assume that $\X \subset \R^n$ is a compact and convex set with the \abc{Bregman} diameter $D_{\X, V}$, function $g$ is convex and $L$-smooth w.r.t. norm $\| \cdot \|$ on $\X$, $f$ is a convex, $G$-Lipschitz w.r.t. norm $\| \cdot \|_2$ on $\X$ and \ab{generally non-differentiable} function. Assume we can use the first-order oracle \ab{(gradients)} for $g(x)$ and the zero-order oracle with in an unbiased stochastic noise for $f(x)$, i.e. we have access to
\begin{equation}
    \label{tilde_f} 
    {f}(x, \xi) \eqdef f(x) + \xi\ab{,}
\end{equation}
where $\xi$ is generated randomly \abc{independently} of the point $x$. Additionally, we assume that the noise is unbiased and bounded:
\begin{equation}
    \label{noise} 
    \E \xi = 0, \quad \E[\xi^2] \leq \sigma^2.
\end{equation}


\subsection{From first to zero-order}

Before presenting the main algorithm, let us understand the properties of the approximation \eqref{grad_f} that we use. Most of these properties have already been encountered in the literature \cite{shamir2017optimal}, we modify only a few of them for our case. 
These properties are associated with the following object
\begin{equation}
    \label{F} 
    F(x) \eqdef \E [f(x + r \tilde e)],
\end{equation}
\abc{where $\tilde e \in \R^n$ is a vector randomly uniformly generated from the unit Euclidean ball and $r > 0$ is the smoothing parameter.}
The function $F$ is called the "smoothed" version of the function $f$. Our algorithm does not use it in any way, but it will be used in the theoretical analysis. 

\abc{
It is worth noting that the use of $F$ as well as the approximation \eqref{grad_f} implies that the function $f$ \abcd{has} to be defined not only on the set $\X$ itself, but also on some neighborhood of it, in order to make the calculation of $f(x+r\tilde e)$ valid for any $x \in \X$ and $\tilde e$ from the unit Euclidean ball. We will not focus much on this fact, which is standard in the literature on gradient-free optimization. Let us only note that there are two ways to fulfill this assumption if the function $f$ is defined strictly on the set $\X$. The first option is to reduce the initial set in the problem \eqref{problem_orig} and consider $\X_{new}$ such that $\X_{new} + B^d_2(r) \subseteq \X$, where $B^d_2(r)$ is the Euclidean ball with a radius of $r$ (see Section 3.2 from \cite{aleks2020gradientfree} for more details). In addition, it is possible to predefine the function $f$ over the entire space $\mathbb{R}^n$. In particular, it can be done as follows: $f_{new}(x) \eqdef \min_{z \in \X} \left[ f(z) + G \| x-z\|\right]$. In Lemma \ref{lem:fnew} from Appendix, we prove that $f_{new}$ is equivalent to $f$ on $\X$, $G$-Lipschitz w.r.t. norm $\| \cdot \|_2$ on $\R^n$ and convex on $\R^n$. 
}

Now let us move to the important technical results. 
\begin{lemma} [see Lemmas 1 and 2 from \cite{beznosikov2019derivative}] \label{lem1}
    $F(x)$ from \eqref{F} is convex, differentiable and it holds that
    \begin{align}
    \sup_{x \in \X} |F(x) - f(x)| &\leq rG \label{orig_sm_main},\\    
    \label{apr_grad_diff_app}
    \E [\tilde f'_r(x, \xi^{\pm})] - \nabla F(x) &= 0, \\
    \label{apr_grad_bounds_app}
    \E [\|\tilde f'_r(x, \xi^{\pm}) \|^2_*]  &\leq p^2(n) \left(8nG^2 + \frac{2n^2\sigma^2}{r^2} \right), 
\end{align}
where $p^2(n) \eqdef \min\{2q - 1, 32 \log n - 8\} n^{\frac{2}{q} - 1}$.
\end{lemma}
Let us discuss these facts. The property \eqref{apr_grad_diff_app} \abcd{means} that the approximation \eqref{grad_f} is an unbiased estimate of the gradient, not of the original function $f$, but of the smoothed function $F$. This means that we can replace the original problem \eqref{problem_orig} with $\min [F(x) + g(x)]$ and now consider the oracle \eqref{grad_f} for $F$ as an unbiased stochastic gradient with a second moment equal to \eqref{apr_grad_bounds_app}. \abcd{The question arises, how much the new problem is different from the original one?} \eqref{orig_sm_main} says that for a small parameter $r$ the original problem \eqref{problem_orig} and the new one are very close. The proof of convergence of the algorithm \abcd{will be built} on this idea.

\subsection{Algorithm and convergence analysis}

As mentioned above, our algorithm is based on the sliding algorithm \cite{lan2016gradient,beznosikov2019derivative}. Our method is a modification of the first-order sliding with a zero-order oracle. The sliding (complexities splitting) effect is achieved by the fact that the method consists of outer and inner loops. At the outer iterations, we compute the gradient of the function $g$, while in the inner loop (prox-sliding procedure), only the function $f$ is used, with fixed information about the gradient \abc{of} $g$.

\begin{algorithm}[h!]
	\caption{One-Point Zero-Order Sliding Algorithm ({\tt OPZOSA})}
	\label{alg}
	\begin{algorithmic}[1]
		\State \textbf{Input:} initial point $x_0 \in \X$, iteration limit $N$ 
            \Statex Let $\{\beta_k\}  > 0$, $\{\gamma_k\}  > 0$, and $\{T_k\} \in \mathbb{N}$ be given
        \Statex Set $\overline x_0 = x_0$
		\For {$k = 1, \ldots, N$}
		\State Set $\underline x_k = (1 - \gamma_k) \overline x_{k-1} + \gamma_k x_{k-1}$ 
        \State \abc{Let $h_k(y) = g(\underline x_k) + \left< \nabla g(\underline x_k), y-\underline x_k\right>$} \label{under_bar_x_k}
		\State Set $(x_k, \Tilde x_k) = \text{{\tt PS}}(h_k, x_{k-1}, \beta_k, T_k)$
		\State Set $\overline x_k = (1 - \gamma_k)\overline x_{k-1} + \gamma_k \Tilde x_k$ \label{bar_x_k}
		\EndFor \\
		\textbf{Output:} $\overline x_N$
	\end{algorithmic}
	{\tt PS} procedure: $(x^+, \Tilde x^+) = \text{{\tt PS}}(h, x, \beta, T)$
	\begin{algorithmic}[1]
		\State \textbf{Input:} function $h: \R^n \to \R$, initial point $x \in \X$, $\beta > 0$, iteration limit $T$
            \Statex Let $\{p_t\} > 0$, $\{\theta_t\} > 0$ be given
            \Statex Set $u_0 = \Tilde u_0 = x$
		\For {$t = 1, \ldots, T$}
		\State $u_t = \argmin\limits_{u \in \X} \left\{ h(u) + \left< \Tilde{f'_r}(u_{t-1}, \xi^{\pm}_{t-1}), u \right> + \beta V(x, u) + \beta p_t V(u_{t-1}, u) \right\}$ \label{u_t}
		\State $\Tilde u_t = (1-\theta_t) \Tilde u_{t-1} + \theta_t u_t$ \label{tilde_u_t}
		\EndFor
		\State \textbf{Output:} $x^{+} = u_T$ and $\Tilde x^{+} = \Tilde u_T$
	\end{algorithmic}
\end{algorithm}

The following theorem gives an estimate for the convergence \abc{rate} of this method:

\begin{theorem}\label{cor:main} Suppose that $ p_t = \frac{t}{2}$, $\theta_t = \frac{2(t+1)}{t(t+3)}$, $\beta_k = \frac{2L}{k}$, $\gamma_k = \frac{2}{k+1}$, 
$$
T_k = \max \left\{ 1; \frac{16Nk^2}{3D_{\X,V}^2 L^2} \cdot \left( 14p^2(n) nG + \frac{p^2(n)n^2\sigma^2}{r^2}\right)\right\}
$$
for $t \geq 1, k \geq 1$. Then for any number of iterations $N$ of Algorithm \ref{alg} it holds that
    \begin{align}
    \label{original_main_conver}
        \mathbb{E}[\Psi_0(\overline x_N)- \Psi_0(x^*)] &\leq 2rG + \frac{20LD_{\X,V}^2}{N(N+1)}.
    \end{align}
Additionally, the total number of {\tt PS} procedure iterations  is 
    \begin{align}
    \label{original_total_main}
        T^{\text{total}} = \frac{(N+1)^4}{D_{\X,V}^2 L^2} \cdot \left( 10p^2(n) nG + \frac{2p^2(n) n^2\sigma^2}{r^2}\right) + N,
    \end{align}
where $p^2(n) \eqdef \min\{2q - 1, 32 \log n - 8\} n^{\frac{2}{q} - 1}$.
\end{theorem}

This theorem shows the significance of the choice of $r$. In particular, it follows from \eqref{original_main_conver} that $r$ should be taken as small as possible. On the other hand, it follows from \eqref{original_total_main} that as $r$ decreases, the total number of internal iterations increases. From here we get a game to some extent: the parameter $r$ must be controlled and adjusted carefully.

\begin{corollary}
Under the assumptions of Theorem \ref{cor:main} and 
if we put $r = \Theta\left(\frac{\varepsilon}{G}\right)$,
then the number of calls of $\nabla g$ and $f$ required by Algorithm 1 to find an $\varepsilon$-solution $\overline x_N$ of \eqref{problem_orig} (i.e. $\E[\Psi_0(\overline x_N)] - \Psi_0(x^*) \le \varepsilon$) is respectively bounded by
\begin{equation*}
    \label{bound_out_main}
    O\left(\sqrt{\frac{L D_{\X,V}^2}{\varepsilon}} \right) \quad \text{and}
    \end{equation*}
\begin{equation*}
    \label{bound_in_main}
  O\left(\sqrt{\frac{L D_{\X,V}^2}{\varepsilon}} + \frac{n \cdot p^2(n) G^2 D_{\X,V}^2}{\varepsilon^2} + \frac{n^2 \cdot p^2(n) G^2\sigma^2 D_{\X,V}^2}{\varepsilon^4} \right),
    \end{equation*}
where $p^2(n) \eqdef \min\{2q - 1, 32 \log n - 8\} n^{\frac{2}{q} - 1}$.
\end{corollary}
\abc{This is the result that we wanted to achieve by the use of the sliding.} \ab{The number of $\nabla g$ calculations}  is not affected in any way by the fact that we use a "very bad" oracle for $f$. Our results also match the obtained bounds for one-point feedback in the non-distributed composite-free case \cite{gasnikov2017stochastic}.

\textbf{Remark.} Note that the second estimate depends on the "geometry" of the problem. In particular, in the Euclidean case $\| \cdot \| = \| \cdot \|_2$ with $q = 2$, we have the following oracle complexity \ab{(number of calls)} for $f$
\begin{equation}
    \label{ecl}
  O\left(\sqrt{\frac{L D_{\X}^2}{\varepsilon}} + \frac{n G^2 D_{\X}^2}{\varepsilon^2} + \frac{n^2 G^2\sigma^2 D_{\X}^2}{\varepsilon^4} \right),
\end{equation}
\ab{where we use $D_{\X}$ as the Euclidean diameter of the set $\X$.}
\abc{A more interesting case is the case when we work in the non-Euclidean setting.  In particular, let us consider a probability simplex with $\| \cdot \| = \| \cdot \|_1$, $q = \infty$ and $D_{\X,V}^2 = 2 \ln n$.} In this case, the estimate is transformed into
\begin{equation}
    \label{necl}
  O\left(\sqrt{\frac{\ln n \cdot L}{\varepsilon}} + \frac{ \ln^2 n \cdot G^2}{\varepsilon^2} + \frac{n \ln^2 n \cdot G^2\sigma^2}{\varepsilon^4} \right).
    \end{equation}
It can be seen that \eqref{necl} improves the estimate \eqref{ecl} $n$ times by using a different geometric setup. Moreover, if the noise $\sigma =0$, our estimates are the same (up to $\ln n$) as the estimates for the full-gradient method \cite{lan2016gradient}.

\subsection{Applications to distributed optimization}

Let us now look at some examples, including those for which sliding gives the estimates necessary in practice. We consider the problems \eqref{centr} and \eqref{decentr}.
\abcd{In the introduction, we briefly mentioned  that in these problems,} we need to reduce the number of $\nabla g$ calls, and thus the number of communications. Indeed, in order to calculate the gradient $g$ in the problem \eqref{centr}, we need to know the value of $\bar x$. \abcd{However, it is impossible to calculate} this using local computations only, it means \abcd{that one needs} help from the central server: send all current $x_m$ and get the average $\bar x$. At the same time, all calculations of $f$ do not require any communication. To compute $f(x_1,\ldots,x_M)$ we need to compute $f_m(x_m)$, and these are just the values of local functions on local variables. For the problem \eqref{decentr}, the same reasoning is valid, but communication takes place with neighbors using the gossip protocol with $W$. \ab{In more details, the computation of $\nabla g$ requires the calculation of the product $WX$.  Since the matrix $W$ has only non-zero weights, when two workers are connected, the multiplication of $WX$ corresponds to the exchange of information between neighbours, where each of the workers, having received parcels from neighbours, averages them with weights according to the matrix $W$. Therefore, the number of $WX$ multiplications is the total number of required communications in the case of a decentralized protocol. Computing $\nabla g$ requires such a multiplication, then the number of communications is equivalent to the calls of the first-order oracle for $g$.}

\abc{
Now we are ready to obtain estimates for the problems \eqref{centr} and \eqref{decentr} from the general results of the previous section. We consider the Euclidean case. It remains only to describe the properties of the problems \eqref{centr} and \eqref{decentr}. For this purpose we give the following lemma. 
\begin{lemma} \label{lem:prop}
Let the functions $f_m$ from \eqref{centr} and \eqref{decentr}  be $G$-Lipschitz w.r.t. $\| \cdot \|_2$, the local noises $\xi_m$ for all $f_m(x_m, \xi_m) = f_m(x_m) + \xi_m$ be independent, unbiased and bounded: $\E \xi_m = 0$, $\E[\xi^2_m] \leq \sigma^2$, the Euclidean diameter of the set $\X$ be equal to $D_{\X}$. Then the following facts for \eqref{centr} and \eqref{decentr} are valid:
\begin{itemize}
    \item $f(x_1,\ldots,x_M)$ is $(G/\sqrt{M})$-Lipschitz w.r.t. $\| \cdot \|_2$ of $\textbf{x} = [x_1^T, \ldots x_M^T]^T$;
    \item $g(x_1, \ldots, x_M)$ from \eqref{centr} is $(\lambda/M)$-smooth w.r.t. $\| \cdot \|_2$ of $\textbf{x}$, and $g$ from \eqref{decentr} is $(\lambda \lambda_{\max}(W)/M)$-smooth w.r.t. $\| \cdot \|_2$ of $\textbf{x}$;
    \item  the Euclidean diameter $D_{\X^M}$ of the set $\X^M$ is equal to $\sqrt{M} D_{\X}$;
    \item the noise of the function $f(x_1,\ldots,x_M, \xi_1, \ldots, \xi_M)$ is unbiased and bounded by $\sigma^2/M$.
\end{itemize}
\end{lemma}
One can note that if $\X \in \R^n$, then the dimension of vector $\textbf{x} = [x_1^T, \ldots x_M^T]^T$ is $n \cdot M$. Recall that the number of computations of $\nabla g(x_1,\ldots,x_M)$ corresponds to the number of communication rounds, and the calls $f(x_1,\ldots,x_M)$ -- to the local gradient-free calculations. Then the following estimates are valid for the number of communications and local iterations to find $\varepsilon$-solution in terms of the target functions:
}

\begin{itemize}
    \item  in the centralized case
\begin{equation*}
    O\left(\sqrt{\frac{\lambda D_{\X}^2}{\varepsilon}} \right) \quad \text{communication rounds} \quad \text{and}
    \end{equation*}
\begin{equation*}
  O\left(\sqrt{\frac{\lambda D_{\X}^2}{\varepsilon}} + \frac{M n G^2 D_{\X}^2}{\varepsilon^2} + \frac{M n^2 G^2\sigma^2 D_{\X}^2}{\varepsilon^4} \right) \quad \text{local computations;}
    \end{equation*}
    \item and in the decentralized case
    \begin{equation*}
    O\left(\sqrt{\frac{\lambda \lambda_{\max}(W) D_{\X}^2}{\varepsilon}} \right) \quad \text{communication rounds} \quad \text{and}
    \end{equation*}
\begin{equation*}
  O\left(\sqrt{\frac{\lambda \lambda_{\max} (W) D_{\X}^2}{\varepsilon}} + \frac{M n G^2 D_{\X}^2}{\varepsilon^2} + \frac{M n^2 G^2\sigma^2 D_{\X}^2}{\varepsilon^4} \right) \quad \text{local computations.}
    \end{equation*}
\end{itemize}
\ab{This is a rather remarkable result. We have "very bad" local functions (they are non-smooth and only with zero-order information), but this fact does not dramatically affect the communication \abc{complexity}. 

\abc{It is important to note that using \eqref{grad_f} for $f(x_1,\ldots,x_M)$ implies that we use a common vector $e$ for the full vector $\textbf{x}$. It turns out that each device for its $x_m$ have to take the corresponding part from the common $e$. This can be accomplished by using the same random generator with the same seed on all devices. In fact, similar results could be obtained using individual and independent directions on all devices. }

It only remains to discuss the choice of \abc{the parameter $\lambda$} for the problems \eqref{centr} and \eqref{decentr}. In fact, this is a key point in personalized learning: a small parameter $\lambda$ is a small contribution of the regularizer to the whole problem, in turn it gives a small penalty for the fact that all the local variables are not similar to each other, it means that each worker takes a little from the other users and relies mostly on the local information. The opposite situation is observed with a large $\lambda$: all $\{x_m\}$ tend to the same value. In particular, there are two extreme cases:}
 \begin{itemize}
    \item If $\lambda = 0$, then \eqref{centr} and \eqref{decentr} become 
    \begin{equation*}
        \min_{X \in \mathcal{X}^M}\sum^M_{m = 1} f_m(x_m).
    \end{equation*}
    This is equivalent to independent local optimization without communications.
   \item
   As $\lambda \to +\infty$,  \eqref{centr} and \eqref{decentr} tend to the distributed problem with equal local arguments:
    \begin{equation*}
        \min_{x_1 = \dots = x_M \in \mathcal{X}}\sum^M_{m = 1} f_m(x_m).
    \end{equation*}
    Note that the infinite $\lambda$ can mess up communications bounds. To solve this issue, $\lambda$ can be taken large but not infinite \cite{dvinskikh2021decentralized,beznosikov2019derivative}. In particular, it is enough to take $\lambda = \frac{G^2}{\lambda^+_{\min} (W) \varepsilon}$, where $\lambda^+_{\min}(W)$ is the smallest positive eigenvalue of the matrix $W$. And then we have the following communication complexities:
    $$
    O\left(\frac{G D_{\X,V}}{\varepsilon} \right), \quad O\left(\sqrt{\frac{\lambda_{\max}(W)}{\lambda^+_{\min}(W)}} \cdot \frac{G D_{\X,V}}{\varepsilon} \right)
    $$
    in the centralized and decentralized cases, respectively.
\end{itemize}

\section{Experiments}

The purpose of our experiments is to compare how our method works in practice in comparison with classical methods. In particular, we compare our method with methods that do not take into account the composite structure of the problem \eqref{problem_orig}. As such a method, we consider Mirror Descent in two settings: in the first case, we consider full-gradient Mirror Descent \cite{nemirovsky1983problem}, which uses \ab{$f' + \nabla g$ as a \abc{subgradient}}; we also consider \abcd{a gradient-free version}, which uses $\tilde f'_r(x,\xi^{\pm}) + \nabla g$ as a "subgradient".

Comparison is made on the problem of distributed computation of geometric median \cite{minsker2015geometric,beznosikov2019derivative}. We have $N$ vectors $\{b_i\}_{i=1}^N \in \R^n$:
\begin{equation}
    \min\limits_{x\in\R^n} f(x) = \sum\limits_{i=1}^N \|x - b_i\|_2.
    \label{eq:geom_median_problem}
\end{equation}
We distributed the vectors $\{b_i\}$ among 10 computing devices. Then the problem can be written in the form \eqref{decentr}:
\begin{equation}
    \min\limits_{X\in\R^{n \times M}}  \underbrace{\sum\limits_{m=1}^M \overbrace{\sum\limits_{i=1}^{\lfloor N/M \rfloor}\|x_m - \abc{\hat b^i_m}\|_2}^{f_m(x_m)}}_{f(X)} + \underbrace{\frac{\lambda}{2}\|\sqrt{W} X\|_2^2}_{g(X)},\label{eq:geom_median_problem_distrib}
\end{equation}
\abc{where $\hat b^i_m = b_{(m-1) \cdot \lfloor N/M \rfloor + i}$.}
\ab{It is easy to verify that each $f_m$ is non-smooth, but $1$-Lipschitz and convex, also, as noted before, $g$ is $\lambda \lambda_{\max}(W)$-smooth and convex. To make our setting stochastic, each time when we call $f'$ or value $f$, we independently generate normal noise vectors $\{\xi_i\}$ for each vector $\{b_i\}$ and compute $f'$ or $f$ using $\{b_i +\xi_i\}$ instead of \abcd{true} $\{b_i\}$.}

\begin{figure*}[h!]
\centering
\includegraphics[width =  0.49\textwidth]{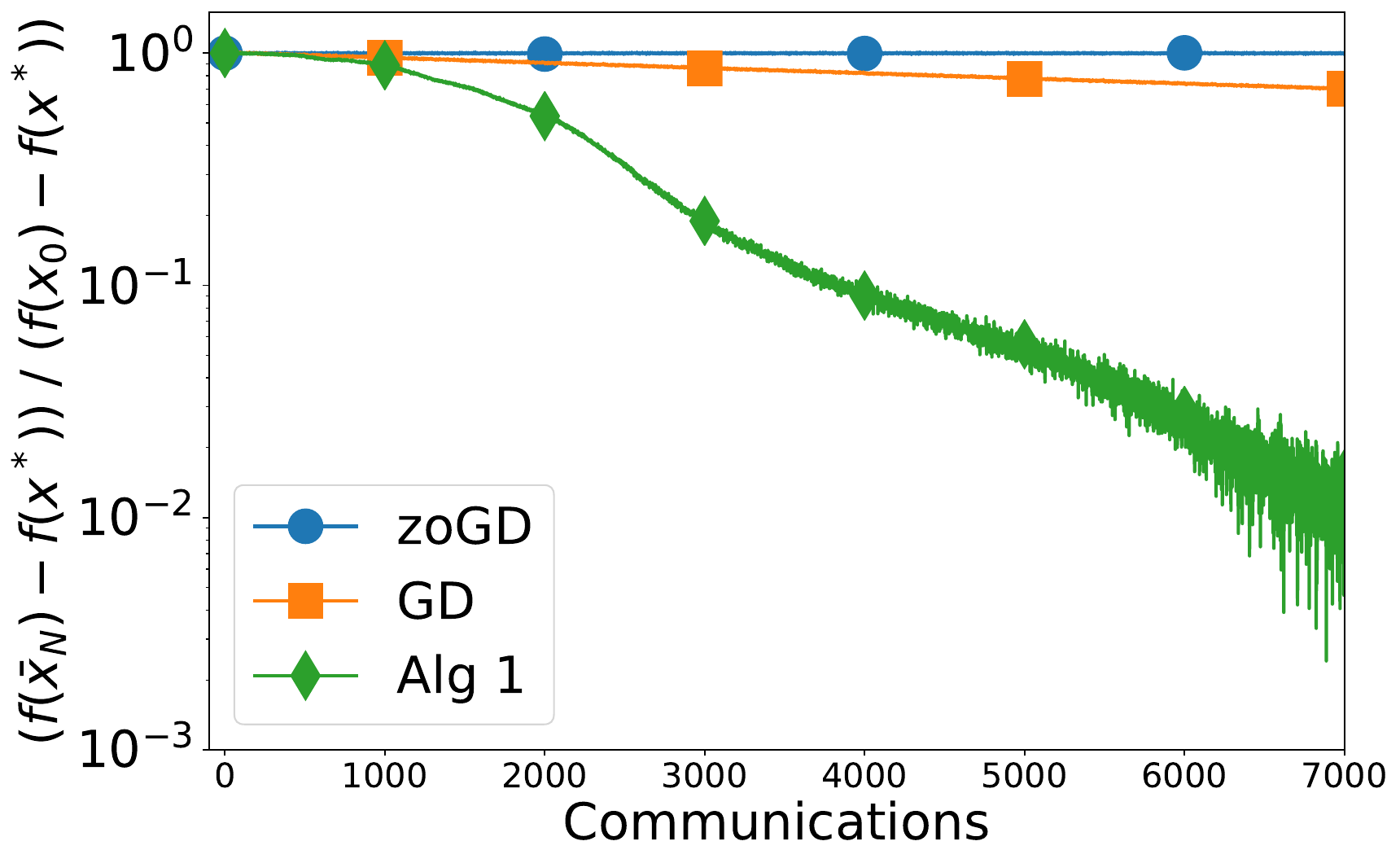}
\includegraphics[width =  0.49\textwidth]{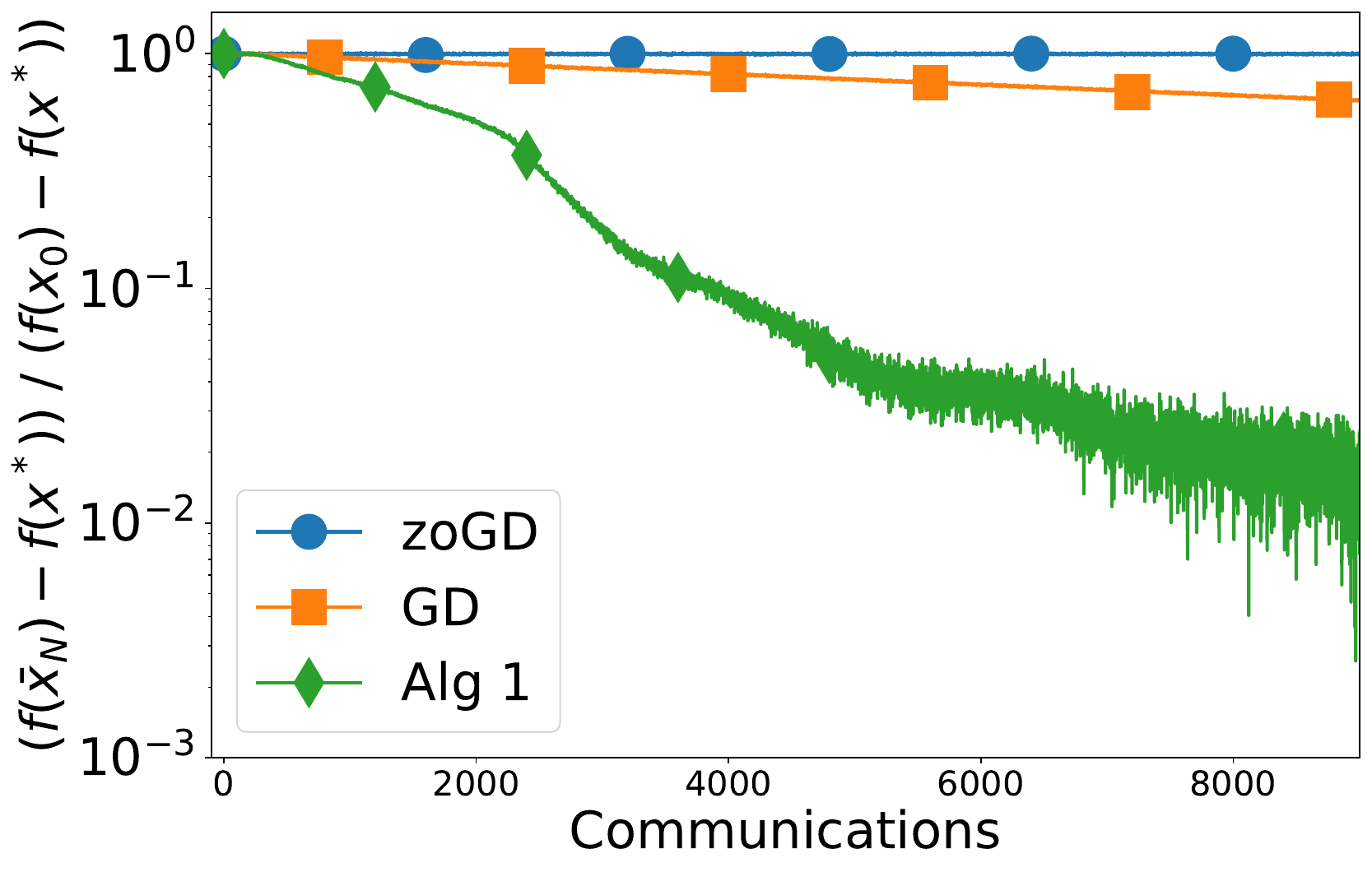}\\
\begin{minipage}{0.50\textwidth}
\centering
~~~~~~(a) star
\end{minipage}%
\begin{minipage}{0.50\textwidth}
\centering
~~~~~~~~(b) complete
\end{minipage}%
\\
\includegraphics[width =  0.49\textwidth]{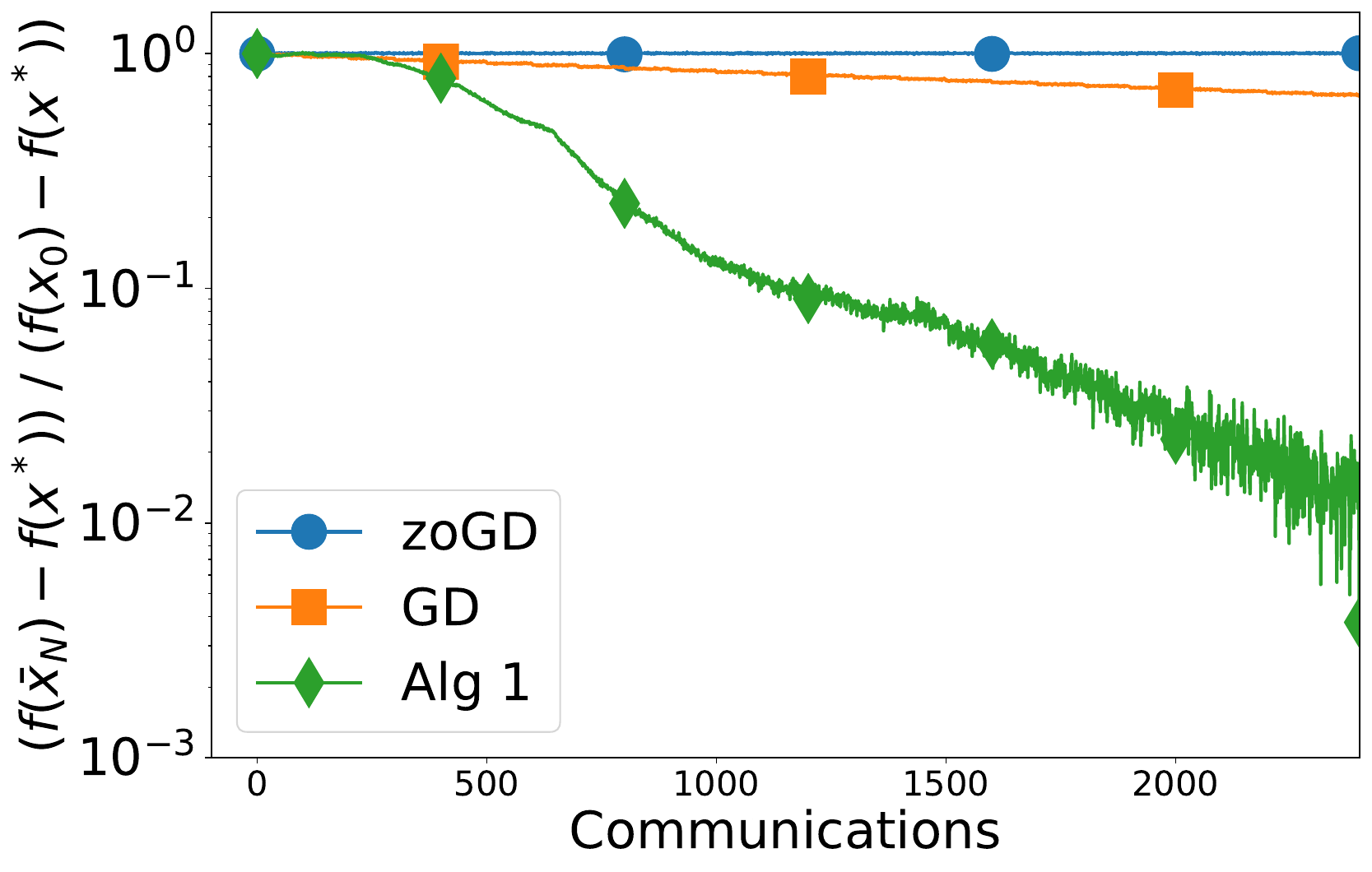}
\includegraphics[width =  0.49\textwidth]{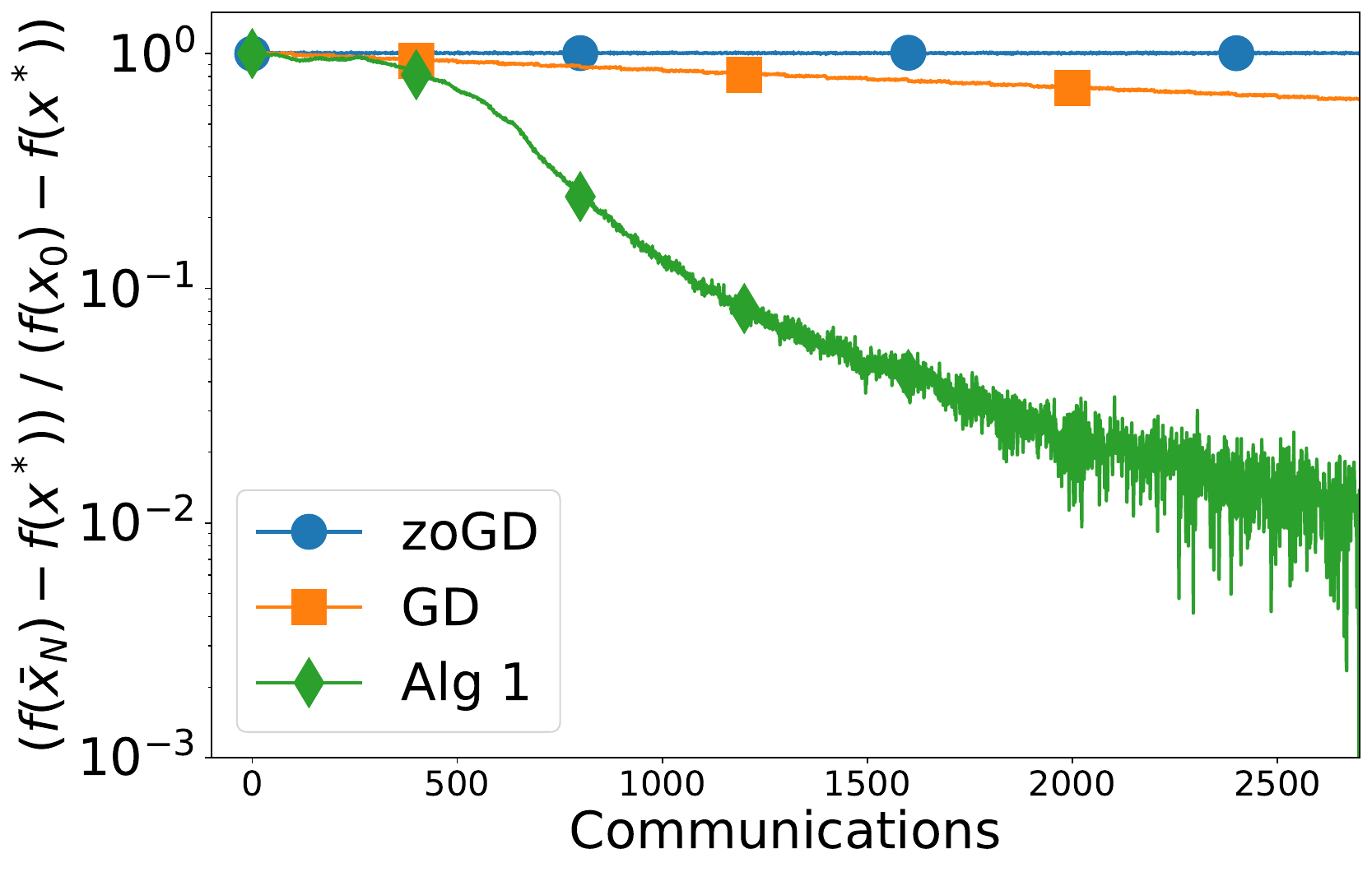}
\begin{minipage}{0.49\textwidth}
\centering
~~~~~(c) chain
\end{minipage}%
\begin{minipage}{0.49\textwidth}
\centering
~~~~~~~(d) cycle
\end{minipage}%
\caption{Algorithm \ref{alg} (green), first-order Mirror Descent (orange) and zero-order Mirror Descent (blue) applied to solve \eqref{eq:geom_median_problem} in distributed form \eqref{eq:geom_median_problem_distrib} with different network topologies. }
\label{fig:distrib_reg10^2}
\end{figure*}

We run \abcd{Algorithm \ref{alg}}, the first-order Mirror Descent  \cite{nemirovsky1983problem}
and the zero-order Mirror Descent \cite{duchi2015optimal} on the problem \eqref{eq:geom_median_problem_distrib} with $n=100$, $N=50$ and $\lambda = 10^2$. 
The vectors $b_1,\ldots, b_N$ are generated as i.i.d. samples from the normal distribution $\cN(\one, 2I_n)$, the noise $\xi$ is also generated from the normal distribution $\cN(0, \sigma^2 I_n)$, where $\sigma = 0.01$. We consider different decentralized  topologies: \abcd{star, complete graph, chain, and cycle}. Both competitor-methods are tuned for better convergence: in the first-order Mirror Descent we tune the step size, in the zero-order Mirror Descent we put $r = 10^{-2}$ in \eqref{grad_f} and also tune the step size. Algorithm \ref{alg} is not tuned and is used with a theoretical set of parameters from Theorem \ref{cor:main} and $r = 10^{-2}$ in \eqref{grad_f}.
The main comparison criterion is the number of communications, i.e. calls to the oracle $\nabla g = \lambda WX$. \ab{We measure the quality of the solution with respect to the target function $f$ of the original problem \eqref{eq:geom_median_problem}, the problem \eqref{eq:geom_median_problem_distrib} is auxiliary and arises due to the distributional nature of the formulation.}  See Figure~\ref{fig:distrib_reg10^2} for the results. We notice that in these \abcd{experiments Algorithm \ref{alg}} outperforms even Mirror Descent which is the first-order method. This shows the importance of taking into account the composite structure of the problem.

\section*{Acknowledgments}

The work on the new version of the paper was  supported by Russian Science Foundation (project No. 23-11-00229).

\bibliographystyle{tfs}
\bibliography{ltr}

\clearpage
\appendix

\part*{Supplementary Material}

\section{Basic Facts}

\begin{lemma}
For an arbitrary integer $m\ge 1$ and an arbitrary set of non-negative numbers $\{a_i\}_{i=1}^m$ the following inequality holds
\begin{equation}
    \left(\sum\limits_{i=1}^m a_i\right)^2 \le m\sum\limits_{i=1}^m a_i^2. \label{eq:squared_sum}
\end{equation}
\end{lemma}

\begin{lemma}[H{\"o}lder's inequality]
For an arbitrary $x,y\in\R^n$ the following inequality holds
\begin{equation}
    \la x, y \ra \le \|x\|_*\cdot\|y\|.
    \label{eq:holder_ineq}
\end{equation}
\end{lemma}

\begin{lemma}[Jensen's inequality]
For an arbitrary convex on $\X \subset \R^n$ function $q$,  an arbitrary integer $m\ge 1$, an arbitrary set of $\{x_i\}_{i=1}^m \in \X$ and an arbitrary set of non-negative numbers $\{a_i\}_{i=1}^m$ such that $\sum_{i=1}^m a_i = 1$ the following inequality holds
\begin{equation}
    q\left( \sum\limits_{i=1}^m a_i x_i\right) \leq \sum\limits_{i=1}^m a_i q\left( x_i\right).
    \label{eq:jensen_ineq}
\end{equation}
\end{lemma}

\begin{lemma}[Cauchy-Schwarz inequality for random variables]
For an arbitrary real valued random variables $\xi$ and $\eta$ such that $\EE[\xi^2] < \infty$ and $\EE[\eta^2] < \infty$ the following inequality holds
\begin{equation}
    \EE[\xi\eta] \le \sqrt{\EE[\xi^2]\EE[\eta^2]}\label{eq:cauchy_schwarz_random}.
\end{equation}
\end{lemma}

\begin{lemma}[Strong convexity of Bregman divergence] \label{eq:bregman_key_property}
For an arbitrary points $x, y \in \X$ the following inequality holds
\begin{equation}
    V(x,y) \ge \frac{1}{2}\|x-y\|^2.
\end{equation}
\end{lemma}

\section{Auxiliary Results}

\begin{lemma}[Lemma 9 from \cite{Shamir15}]\label{lem:lemma_9_shamir} Let $q$ be a $Q$-Lipschitz \abcd{w.r.t.} norm $\| \cdot \|_2$ function and $e$ be uniformly distributed on the Euclidean unit sphere, then 
\begin{equation*}
    \sqrt{\mathbb{E}[(q(e) - \mathbb{E}q(e))^4]} \leq \frac{4 Q^2}{n}.
\end{equation*}
\end{lemma}

\begin{lemma}[Lemma 3.5 from \cite{lan}]\label{lem:lemma3.5_lan} Let a convex function $q: \X \to \mathbb{R}$, points $\tilde x, \tilde y \in \X$ and scalars $\mu_1, \mu_2 \geq 0$ be given. Let $\nu: \X \to \mathbb{R}$ be a differentiable convex function and 
$
V(x,z) = \nu(z) - [\nu(x) + \abc{\langle \nabla \nu(x), z-x \rangle}].
$
If
$
\tilde u = \argmin_{u \in \X}\{q(u) + \mu_1 V(\tilde x, u) + \mu_2 V(\tilde y, u)\},
$
then for any $u \in \X$ we have
\begin{equation*}
q(\tilde u) + \mu_1 V(\tilde x, \tilde u) + \mu_2 V(\tilde y, \tilde u) \leq q(u) + \mu_1 V(\tilde x,  u) + \mu_2 V(\tilde y,  u) - (\mu_1 + \mu_2)V(\tilde u, u).
\end{equation*}
\end{lemma}

\begin{lemma}[Lemma 3.17 from \cite{lan}]\label{lem:lemma3.17_lan}
 Let $w_k \in (0;1)$, $k \geq 1$ be given. Also let us denote
\begin{align} W_k = 
\begin{cases} 
    1, & k=1,\\
    (1 - w_k)W_{k-1}, & k>1.
\end{cases}\nonumber
\end{align}
Suppose that the sequence $\{\Delta_k\}_{k \geq 0}$ satisfies
$$
\Delta_k \leq (1 - w_k)\Delta_{k-1} + B_k
$$
for all $k \geq 1$ and some positive constants $\{B_k\}_{k  \geq 1}$. Then for all $k \in \mathbb{N}$ we have 
\begin{equation*}
\Delta_k \leq W_k
(1 - w_1)
\Delta_{0} + W_k \sum\limits_{i=1}^{k} \frac{B_i}{W_i} \leq W_k
\Delta_{0} + W_k \sum\limits_{i=1}^{k} \frac{B_i}{W_i}.
\end{equation*}
\end{lemma}

\begin{lemma}\label{lem:bounded_grad_to_lan's_condition}
    Assume that for the differentiable function $q$ defined on a closed and convex set $\X$ there exists $Q$ such that for all $x \in \X$
    \begin{eqnarray}
    \label{bound_grad} 
     \|\nabla q(x) \|_* \leq Q.
\end{eqnarray}
Then it holds that for all $x,y \in \X$
\begin{equation*}
    q(x) \le q(y) + \langle \nabla q(y), x - y \rangle + 2Q\|x-y\|.
\end{equation*}
\end{lemma}
\textit{Proof of Lemma~\ref{lem:bounded_grad_to_lan's_condition}:}
    For an arbitrary points $x,y\in \X$ we have
    \begin{align*}
    q(x) &= q(y) + \int\limits_{0}^{1} \langle \nabla q(y + \tau(x-y)), x-y \rangle d\tau 
     \nonumber\\
     &= q(y) + \langle \nabla q(y), x-y \rangle + \int\limits_0^1 \langle \nabla q(y + \tau(x-y)) - \nabla q(y), x-y \rangle d\tau.
    \end{align*}
Using \eqref{eq:holder_ineq} and then \eqref{bound_grad}, we obtain 
    \begin{align*}
    q(x) &\overset{}{\leq} q(y) + \langle \nabla q(y), x-y \rangle + \int\limits_0^1 \|\nabla q(y + \tau(x-y)) - \nabla q(y) \|_* \cdot \|x-y\| d\tau 
    \nonumber\\
    &\overset{}{\leq} q(y) + \langle \nabla q(y), x-y \rangle + \int\limits_0^1 2Q\|x-y\| d\tau 
    \nonumber\\
    &\leq q(y) + \langle \nabla q(y), x - y \rangle + 2Q\|x-y\|.
\end{align*}

\abc{
\section{Lemma on $f_{new}$}

\begin{lemma} \label{lem:fnew} 
Let $f_{new}$ be defined as follows
    \begin{align*}
        f_{new}(x) = \min_{z \in \X} \left[ f(z) + G \| x-z\|\right]
    \end{align*}
Then $f_{new}$ is equivalent to $f$ on $\X$, $G$-Lipschitz w.r.t. norm $\| \cdot \|_2$ on $\R^n$ and convex on $\R^n$.
\end{lemma}

\textit{Proof of Lemma~\ref{lem:fnew}:}
Since $f$ is Lipschitz continuous on the compact set $\X$, it follows that for any $x \in \R^d$ the function 
$
\phi_x(z) \eqdef [f(z) + G \| x-z\|]
$
is continuous on $\X$. Then this function reaches its optimal point on $\X$, i.e. for any $x \in \R^d$ there exists $z_x \in \X$ such that $f_{new}(x) = \phi_x (z_x) = \min_{z \in \X} [f(z) + G \| x-z\|]$. Now we are ready to prove the lemma.

We start from the fact that $f_{new}$ is equivalent to $f$ on $\X$. 
It is easy to see that for all $x \in \X$
$$
f_{new} (x) = \min_{z \in \X}  \phi_x (z) \leq \phi_x (x) = f(x) + G \| x-x\| = f(x).
$$
Let us prove by contradiction and suppose that there exists $\tilde x \in X$ such that $f_{new} (\tilde x) < f(\tilde x)$, i.e. 
\begin{align*}
    f(z_{\tilde x}) +  G \| z_{\tilde x} - \tilde x\| =  \min_{z \in \X} \phi_{\tilde x} (z) = f_{new} (\tilde x) < f(\tilde x).
\end{align*}
Here we also used the discourse above on the existence of $z_{\tilde x}$. 
From this estimate, we have
$$
G \| z_{\tilde x} - \tilde x\| < f(\tilde x) - f(z_{\tilde x}). 
$$
\abcd{On the other hand, the $G$-Lipshitzeness of the function $f$ implies that} 
\begin{align*}
   f(\tilde x) - f(z_{\tilde x}) \leq G \| z_{\tilde x} - \tilde x\|.
\end{align*}
We come to a contradiction. It means that $f_{new} (x) = f(x)$ for all $x \in \X$. 

Next, we prove that $f_{new}$ is $G$-Lipschitz. Let us consider any $x,y \in \R^n$ and without loss of generality assume that $f_{new}(x) \geq f_{new}(y)$. Then, we get
\begin{align*}
|f_{new}(x) - f_{new}(y)|
\leq
f_{new}(x) - f_{new}(y)
=
\phi_x (z_x) - \phi_y (z_y).
\end{align*}
 One can note that $\phi_x (z_x) = \min_{z \in \X} \phi_x (z) \leq \phi_x (z_y)$. Combining these two facts, we have
\begin{align*}
|f_{new}(x) - f_{new}(y)|
&\leq
\phi_x (z_y) - \phi_y (z_y)
\\
&= 
f(z_y) + G \| x - z_y\| -  \left[f(z_y) + G \| y - z_y\|\right]
\\
&= G \left[\| x - z_y\| - \| y - z_y\|\right]
\\
&\leq G \| x - y\|.
\end{align*}
This is what was required.

Finally, we prove the convexity of $f_{new}$ by Definition \ref{def:conv}. For the non-differentiable function $f_{new}$ we consider any $\alpha \in [0;1]$ and get
\begin{align*}
\alpha f_{new}(x) + (1 - &\alpha)f_{new}(y)
\\
&= 
\alpha \left[ f(z_x) + G \| x - z_x\|\right] +  (1 - \alpha)\left[f(z_y) + G \| y - z_y\|\right]
\\
&=
\alpha f(z_x) + (1 - \alpha) f(z_y) + \alpha G \| x - z_x\| + (1 - \alpha) G \| y - z_y\|
\\
&\geq f(\alpha z_x + (1 - \alpha) z_y) + G \| \alpha x + (1 - \alpha) y - \alpha z_x - (1 - \alpha)z_y\|
\\
&\geq 
f_{new} (\alpha x + (1 - \alpha) y).
\end{align*}
In the third step, we used the convexity of $f$ and $\| \cdot\|$. In the last step, we took into account the definition of $f_{new}$: $f(\alpha z_x + (1 - \alpha) z_y) + G \| \alpha x + (1 - \alpha) y - \alpha z_x - (1 - \alpha)z_y\| \geq \min_{z \in \X} [f(z) + G \| \alpha x + (1 - \alpha) y -z\|] = f_{new}(\alpha x + (1 - \alpha) y)$. This finishes the proof of the third fact. 
\EndProof
}

\section{Proof of Lemma \ref{lem1}}

In this section, we prove Lemma \ref{lem1}. For convenience, we divided the proof into two lemmas. Lemma \ref{first_lemma} gives the properties of the function \eqref{F}, and Lemma \ref{second_lemma} -- the properties of the approximation \eqref{grad_f}. Also, for convenience, we duplicate the statements of Lemma \ref{lem1}.

\begin{lemma} [Lemma 8 from \cite{shamir2017optimal} and Lemma 1 from \cite{beznosikov2019derivative}] \label{first_lemma}
    $F(x)$ defined in \eqref{F} is a convex, differentiable function. Moreover, $F(x)$ satisfies
    \begin{align}
        \nabla F(x) &= \mathbb{E} \left[\frac{n}{r} f(x + r e)e\right], \label{grad_sm_main_app}\\
        \sup_{x \in \X} |F(x) - f(x)| &\leq rG \label{orig_sm_main_app},\\    
        \|\nabla F(x)\|_* &\leq 2 \sqrt{n} p(n)  G, \label{eq:norm_nabla_F_bound_main_app}
    \end{align}
    where $p^2(n)$ defined as follows: $\sqrt{\E[\|e\|_*^4]} \leq p^2(n) \eqdef \min\{2q - 1, 32 \log n - 8\} n^{\frac{2}{q} - 1})$.
\end{lemma}
\textit{Proof of Lemma \ref{first_lemma}:} 
The convexity and differentiability of the function $F(x)$ and \eqref{grad_sm_main_app} follows from Lemma 8 of \cite{shamir2017optimal}. Then using sequentially the definition of $F(x)$, the properties of the expectation associated with the absolute value and $G$-Lipschitzness of $f$, we get that for all $x \in \X$
\begin{align*}
    |F(x) - f(x)| &= |\E [f(x + r \tilde e)] - f(x)| \\
    &\leq \E [|f(x + r \tilde e) - f(x)|] \\
    &\leq \E [G \cdot \|r \tilde e\|_2]
    \\
    &\leq Gr.
\end{align*}
It proves \eqref{orig_sm_main_app}. Finally, we deal with \eqref{eq:norm_nabla_F_bound_main_app}. By the symmetry of the distribution of $e$ and  \ab{\eqref{grad_sm_main_app}}, we get:
\begin{align*}
    \|\nabla F(x)\|_*^{2} &= \left\| \E \left[ \frac{n}{r} f(x + re) e\right] \right\|_*^{2} \\ 
    &= \left\| \E \left[ \frac{n}{2r} f(x + re) e \right] - \E \left[\frac{n}{2r} f(x - re) e\right]\right\|_*^{2}\\
    &\leq \frac{n^2}{4r^2} \E \left[(f(x + re) - f(x - re))^2 \|e\|_*^2\right]\\ 
    &= \frac{n^2}{4r^2} \E \left[((f(x + re) - \alpha) - (f(x - re) - \alpha))^2 \|e\|_*^2 \right] ,
\end{align*}
\ab{where $\alpha$ is some constant, which we will define later.}
Next, we apply \eqref{eq:squared_sum} and obtain:
\begin{align*}
    \|\nabla F(x)\|_*^{2} \leq \frac{n^2}{2r^2} \E \left[(f(x + re) - \alpha)^2\|e\|_*^2 \right] + \frac{n^2}{2r^2}\E \left[(f(x - re) - \alpha)^2\|e\|_*^2 \right].
\end{align*}
Since the distribution of $e$ is symmetric, one can note that $\E \left[(f(x + re) - \alpha)^2\|e\|_*^2 \right] = \E \left[(f(x - re) - \alpha)^2\|e\|_*^2 \right]$ and then
\begin{equation*}
    \|\nabla F(x)\|_*^{2} \leq \frac{n^2}{r^2} \E \left[(f(x + re) - \alpha)^2\|e\|_*^2 \right].
\end{equation*}
Using the Cauchy-Schwarz inequality \eqref{eq:cauchy_schwarz_random}, we get
\begin{equation*}
    \|\nabla F(x)\|_*^{2} \leq \frac{n^2}{r^2} \E[(f(x + re) - \alpha)^2\|e\|_*^2] \leq \frac{n^2}{r^2} \sqrt{\E[\|e\|_*^4]\E[(f(x + re) - \alpha)^4]}.
\end{equation*}
With the notation of $p(n)$ from the statement of the lemma, we have
\begin{equation}
    \label{eq:24_2}
    \|\nabla F(x)\|_*^{2} \leq \frac{n^2 p^2(n)}{r^2} \sqrt{\E[(f(x + re) - \alpha)^4]}.
\end{equation}
Taking \abcd{$\alpha = \mathbb{E}[f(x + re)]$}, having that $f(x + re)$ is $Gr$-Lipshitz w.r.t. $e$ in terms of $\|\cdot\|_2$ and using Lemma \ref{lem:lemma_9_shamir}, we get:
\begin{equation}
    \label{eq:24_1}
    \sqrt{\E[(f(x + re) - \alpha)^4]} \leq \frac{4G^2 r^2}{n}.
\end{equation}
Combining \eqref{eq:24_2} and \eqref{eq:24_1}, we prove that
\begin{equation*}
    \|\nabla F(x)\|^2_* \le \ab{4} n p^2(n) G^2.
\end{equation*}
\EndProof

\begin{lemma}[see Lemma 2 from \cite{beznosikov2019derivative}] \label{second_lemma}
For $\tilde f'_r$ given in \eqref{grad_f} the following \abc{relations} hold:
\begin{align}
    \label{apr_grad_diff}
    \E [\tilde f'_r(x, \ab{\xi^{\pm}})] - \nabla F(x) &= 0, \\
    \label{apr_grad_bounds}
    \E [\|\tilde f'_r(x, \ab{\xi^{\pm}}) \|^2_*]  &\leq p^2(n) \left(8nG^2 + \frac{2n^2\sigma^2}{r^2} \right),
\end{align}
where $p^2(n)$ satisfies $ \sqrt{\E[\|e\|^4_*]} \le p^2 (n)$.
\end{lemma}
\textit{Proof of Lemma \ref{second_lemma}:} 
We start from \eqref{apr_grad_diff}. With the definition \eqref{grad_f}, we get
\begin{align*}
        \mathbb{E}[\tilde{f}_r'(x, \xi^{\pm})] &= \frac{n}{2r}\left(\mathbb{E}[f(x + re, \xi^+) e] - \mathbb{E}[f(x - re, \xi^-) e]\right) \\
        &= \frac{n}{2r}\left(\mathbb{E}[f(x + re) e] - \mathbb{E}[f(x - re) e]\right) + \frac{n}{2r}\mathbb{E}[(\xi^+ - \xi^-) e]
\end{align*}
Taking into account the independence of $e$, $\xi$ and \eqref{noise} we have $\mathbb{E}[(\xi^+ - \xi^-) e] = \mathbb{E}_e\left[\mathbb{E}_\xi[(\xi^+ - \xi^-) e]\right] = 0$. Then, using \eqref{grad_sm_main_app}, we obtain
\begin{align*}
        \mathbb{E}[\tilde{f}_r'(x, \xi^{\pm})] - \nabla F(x)
        &= \frac{n}{2r}\mathbb{E}\left[f(x + re) e\right] - \frac{n}{2r} \mathbb{E} \left[f(x - re) e\right] - \nabla F(x) \\
        &= \frac{n}{r}\mathbb{E}\left[f(x + re) e\right] - \nabla F(x) = 0
\end{align*}
Next, we prove \eqref{apr_grad_bounds}. Again with the definition \eqref{grad_f}, we get
\begin{align*}
    \E[ \|\tilde f'_r(x, \xi^{\pm}) \|^2_* ] &=
    \E \left[ \left\| \frac{n}{2r}(f(x + re, \xi^{+}) -  f(x - re, \xi^{-}))e \right\|^2_* \right] \\
    &= \frac{n^2}{4r^2}\E \left[ \left\|(f(x + re) +  \xi^{+} - f(x - re) -  \xi^{-})e \right\| ^2_* \right].
\end{align*}
Using the property \eqref{eq:squared_sum} twice, we get
\begin{align*}
    \E\left[ \|\tilde f'_r(x, \xi^{\pm}) \|^2_* \right] \leq& \frac{n^2}{2r^2}\E\left[ \|(f(x + re) - f(x - re))e\|^2_*\right] +  \frac{n^2}{2r^2}\E\left[ \|(\xi^{+} - \xi^{-})e\|^2_*\right]\\
    \leq& \frac{n^2}{2r^2}\E\left[(f(x + re) - f(x - re))^2\|e\|^2_*\right] \\
    &+  \frac{n^2}{r^2}\E\left[ ((\xi^{+})^2 + (\xi^{-})^2)\|e\|^2_*\right] \\
    =& \frac{n^2}{2r^2}\E\left[(f(x + re) -\alpha - f(x + re) + \alpha)^2\|e\|^2_*\right] \\
    &+  \frac{n^2}{r^2}\E\left[ ((\xi^{+})^2 + (\xi^{-})^2)\|e\|^2_*\right],
\end{align*}
where as in Lemma \ref{first_lemma} we introduce $\alpha$, which we will later set equal to \abcd{$\mathbb{E}_{e}[f(x + re)]$.}
By independence of $\xi^{\pm}$ and $e$, we have
\begin{align*}
    \E\left[ \|\tilde f'_r(x, \xi^{\pm}) \|^2_* \right]
    \leq& \frac{n^2}{2r^2}\E_{e}\left[(f(x + re) -\alpha - f(x + re) + \alpha)^2\|e\|^2_*\right] \\
    &+ \frac{n^2}{r^2}\E_{\xi}\left[\E_{e}\left[ ((\xi^{+})^2 + (\xi^{-})^2)\|e\|^2_* \right]\right] \\
    \leq& \frac{n^2}{r^2}\E_{e}\left[ \left[(f(x + re) -\alpha)^2 + (f(x - re) - \alpha)^2 \right] \|e\|^2_* \right] \\
    &+ \frac{n^2}{r^2}\E_{\xi}\left[\E_{e}\left[ ((\xi^{+})^2 + (\xi^{-})^2)\|e\|^2_* \right]\right].
\end{align*}
Taking into account the symetric distribution of $e$, also using the Cauchy-Schwarz inequality \eqref{eq:cauchy_schwarz_random}, the definition of $p^2(n)$ and \eqref{noise}, we get
\begin{align*}
    \E\left[ \|\tilde f'_r(x, \xi^{\pm}) \|^2_* \right]
    \leq& \frac{\ab{2}n^2}{r^2}\E_{e}\left[(f(x + re) -\alpha)^2\|e\|^2_*\right] + \frac{n^2}{r^2}\E_{\xi}\left[\E_{e}\left[ ((\xi^{+})^2 + (\xi^{-})^2)\|e\|^2_* \right]\right] \\
    \leq& \frac{\ab{2}n^2}{r^2}\sqrt{\E_{e}\left[(f(x - re) - \alpha)^4\right]}\sqrt{\E_{e}\left[\|e\|^4_*\right]} \\
    &+ \frac{n^2}{r^2}\E_{\xi}\left[\E_{e}\left[ ((\xi^{+})^2 + (\xi^{-})^2)\|e\|^2_* \right]\right] \\
    \leq& \frac{\ab{2}n^2p^2(n)}{r^2}\sqrt{\E_{e}[(f(x - re) - \alpha)^4]} + \frac{2n^2p^2(n)\sigma^2}{r^2}.
\end{align*}
Putting \abcd{$\alpha = \mathbb{E}[(f(x + re)]$,} taking into account that $f(x + re)$ is $Gr$-Lipshitz w.r.t. $e$ in terms of $\|\cdot\|_2$ and using Lemma \ref{lem:lemma_9_shamir}, we get:
\begin{equation*}
    \E[ \|\tilde f'_r(x, \xi^{\pm}) \|^2_*]  \le p^2(n)\left(\ab{8}n \ab{G}^2 + \frac{2n^2\sigma^2}{r^2}\right).
\end{equation*}
\EndProof

\section{Proof of Theorem \ref{cor:main}}

In this section, we prove the main theorem. The analysis is based on \cite{lan2016gradient,beznosikov2019derivative}. \abc{Let us consider the following lemma, which provides a way to analyze {\tt PS} \ab{procedure} from Algorithm \ref{alg}.}

\begin{lemma} \label{third_lemma}
Assume that $\{p_t\}_{t\ge 1}$ and $\{\theta_t\}_{t\ge 1}$ in the subroutine {\tt PS} \ab{with the input $(h, x, \beta, T)$}  satisfy
    \begin{align}
        \label{p_t_main} 
        P_t &=
        \begin{cases}
            1 & t = 0, \\
            p_t(1 + p_t)^{-1}P_{t-1} & t \geq 1,
        \end{cases}
        \\
        \label{theta_t_main} 
        \theta_t &= \frac{P_{t-1} - P_t}{(1 - P_t)P_{t-1}} \in [0;1].
    \end{align}
    Then for any $t \geq 1$ and $u \in \X$:
    \begin{align}
        \label{lemma_2_main}
        \frac{\beta}{1 - P_{t}} &V(u_t, u) + (\Phi(\tilde u_t)  - \Phi(u)) \notag\\
    \leq& \frac{\beta P_{t}}{1 - P_{t}} V(x, u) +  \frac{P_{t}}{1 - P_{t}}\sum_{i=1}^t \frac{1}{P_{i-1} p_i} \left[ \frac{(\abc{2} \ab{\tilde G} + \|\delta_i\|_*)^2}{2 \beta p_i } 
    + \langle \delta_i, u - u_{i-1} \rangle \right],
    \end{align}
    where 
    \begin{align}
        \label{Phi_main}
        \Phi(u) &\eqdef h(u) + F(u) + \beta V(x,u), \\
        \label{delta_main}
        \delta_t &\eqdef \tilde f_r' (u_{t-1}, \xi^{\pm}_{t-1}) - \nabla F(u_{t-1}), \\
        \tilde G &\eqdef 2p(n)\sqrt{n}G. \label{tG}
    \end{align}
\end{lemma}
\textit{Proof of Lemma \ref{third_lemma}:}
Lemma \ref{first_lemma} guarantees that $F$ is differentiable and has $\| \nabla F(x)\|_* \leq \tilde G$ from \eqref{tG}. It means we can use Lemma \ref{lem:bounded_grad_to_lan's_condition} for $F$ to get
\begin{align*}
     F(u_t) &\leq F(u_{t-1}) + \langle \nabla F(u_{t-1}), u_t - u_{t-1} \rangle + \abc{2} \tilde G\| u_t - u_{t-1}\|.
\end{align*}
Adding $h(u_t) + \beta V(x, u_t)$ to this inequality and applying \eqref{Phi_main}, we obtain:
\begin{align*}
    \Phi(u_t) 
    =& h(u_t) + F(u_t) + \beta V(x,u_t)\\
    \leq& h(u_t) + F(u_{t-1}) + \langle \nabla F(u_{t-1}), u_t - u_{t-1} \rangle + \beta V(x, u_t) + \abc{2} \tilde G \|u_t - u_{t-1}\|\\
    =& h(u_t) + F(u_{t-1}) + \langle \tilde f'_r(u_{t-1}, \xi^{\pm}_{t-1}), u_t - u_{t-1} \rangle - \langle \delta_t, u_t- u_{t-1} \rangle
    \\
    &+ \beta V(x, u_t) + \abc{2} \tilde G \|u_t - u_{t-1}\|.
\end{align*}
In the last step, we also took into account \eqref{delta_main}, we get
\begin{align}
    \label{temp1}
    \begin{split}
    \Phi(u_t) \le& h(u_t) + F(u_{t-1}) + \langle \tilde f'_r(u_{t-1}, \xi^{\pm}_{t-1}), u_t - u_{t-1} \rangle 
    \\
    &+ \beta V(x, u_t) + (\abc{2} \tilde G + \|\delta_t\|_*) \|u_t - u_{t-1}\|.
    \end{split}
\end{align}
Next, we apply Lemma \ref{lem:lemma3.5_lan} to Line \ref{u_t} of {\tt PS} procedure. Here \abc{$q (\cdot) = h(\cdot) + \langle \Tilde{f'_r}(u_{t-1}, \xi^{\pm}_{t-1}), \cdot \rangle$}, $\mu_1 = \beta$, $\mu_2 = \beta p_t$, $\tilde u = u_t$, $\tilde x = x$ and $\tilde y = u_{t-1}$.
Then we obtain that for all $u \in \X$
\begin{align*}
    h(u_t) + \langle \tilde f'_r(u_{t-1}, & \xi^{\pm}_{t-1}), u_t - u_{t-1} \rangle 
    + \beta V(x, u_t)+ \beta p_t V(u_{t-1}, u_t) 
    \notag\\
    \le& h(u) + \langle \tilde f'_r(u_{t-1}, \xi^{\pm}_{t-1}), u - u_{t-1} \rangle 
    \\
    &+ \beta V(x, u) + \beta p_t V(u_{t-1}, u) - \beta(1 + p_t) V(u_t, u)
    \\
    =& h(u) + \langle \nabla F(u_{t-1}), u - u_{t-1} \rangle + \langle \delta_t, u - u_{t-1} \rangle 
    \\
    &+ \beta V(x, u) + \beta p_t V(u_{t-1}, u) - \beta(1 + p_t) V(u_t, u).
\end{align*}
The convexity of $F$ (see Lemma \ref{first_lemma}) gives that $\langle \nabla F(u_{t-1}), u - u_{t-1} \rangle \leq F(u) - F(u_{t-1})$ and then
\begin{align}
    \label{temp2}
    h(u_t) + F(u_{t-1}) + & \langle \tilde f'_r(u_{t-1}, \xi^{\pm}_{t-1}), u_t - u_{t-1} \rangle 
    + \beta V(x, u_t)+ \beta p_t V(u_{t-1}, u_t) 
    \notag\\ 
    \le& h(u) + F(u) + \langle \delta_t, u - u_{t-1} \rangle + \beta V(x, u) \notag\\
    &+ \beta p_t V(u_{t-1}, u) - \beta(1 + p_t) V(u_t, u)\notag\\
    \leq& \Phi(u) + \beta p_t V(u_{t-1}, u) - \beta (1 + p_t) V(u_t, u) + \langle \delta_t, u - u_{t-1} \rangle.
\end{align}
Here we also used the definition \eqref{Phi_main}. Summing \eqref{temp1} and \eqref{temp2}, one can obtain
\begin{align}
    \label{eq:24_6}
    \begin{split}
    \Phi(u_t) + \beta p_t V(u_{t-1}, u_t) \le& \Phi(u) + \beta p_t V(u_{t-1}, u) - \beta (1 + p_t) V(u_t, u)
    \\
    &+ \langle \delta_t, u - u_{t-1} \rangle  + (\abc{2} \tilde G + \|\delta_t\|_*) \|u_t - u_{t-1}\|.
    \end{split}
\end{align}
Moreover, the strong convexity of $V$ (Lemma \ref{eq:bregman_key_property}) implies that
\begin{align}
    \label{temp3}
    (\abcd{2} \tilde G+\|\delta_t\|_*) \| u_t - u_{t-1}\| \leq& \frac{\beta p_t }{2} \|u_t - u_{t-1}\|^2 + \frac{1}{2 \beta p_t } (\abc{2} \tilde G + \|\delta_t\|_*)^2\notag \\
    \leq& \beta p_t V(u_{t-1}, u_t) + \frac{(\abc{2} \tilde G+\|\delta_t\|_*)^2}{2\beta p_t}.
\end{align}
Combining \eqref{eq:24_6} and \eqref{temp3}, we get
\begin{equation*}
    \Phi(u_t)  - \Phi(u)  \leq  \beta p_t V(u_{t-1}, u) - \beta (1 + p_t) V(u_t, u)
    + \frac{(\abcd{2} \tilde G + \|\delta_t\|_*)^2}{2\beta p_t}
    + \langle \delta_t, u - u_{t-1} \rangle.
\end{equation*}
Now dividing both sides of the above inequality by $1 + p_t$ and rearranging the terms, we get
\begin{equation*}
    \beta V(u_t, u) 
    \leq \frac{\beta p_t}{1+ p_t} V(u_{t-1}, u) + \frac{(\abc{2} \tilde G + \|\delta_t\|_*)^2}{2 \beta (1+p_t) p_t }
    + \frac{\langle \delta_t, u - u_{t-1} \rangle}{1 + p_t} - \frac{\Phi(u_t)  - \Phi(u)} {1+p_t}.
\end{equation*}
Next, we apply Lemma \ref{lem:lemma3.17_lan} with $w_k = 1 -\tfrac{p_t}{1+ p_t} \in (0;1)$, $W_k = P_t$ (see \eqref{p_t_main}) and $\Delta_k = \beta V(u_t, u)$ and get
\begin{align*}
    \beta V(u_t, u) 
    \leq& P_t \beta V(u_{0}, u) +  P_t \sum_{i=1}^t \left[ \frac{(\abc{2} \tilde G + \|\delta_i\|_*)^2}{2 \beta P_i (1+p_i) p_i } 
    + \frac{\langle \delta_i, u - u_{i-1} \rangle}{P_i (1 + p_i)} - \frac{\Phi(u_i)  - \Phi(u)} {P_i (1+p_i)}  \right]. 
\end{align*}
Multiplying by $\tfrac{1}{1 - P_{t}}$ and making rearrangements, we obtain
\begin{align}
    \label{temp21}
    \begin{split}
    \frac{\beta}{1 - P_{t}} &V(u_t, u) + \sum_{i=1}^t \frac{P_{t}}{P_i (1+p_i) (1 - P_{t}) }\cdot (\Phi(u_i)  - \Phi(u)) 
    \\
    \leq& \frac{\beta P_{t}}{1 - P_{t}} V(u_{0}, u) +  \frac{P_{t}}{1 - P_{t}}\sum_{i=1}^t \left[ \frac{(\abc{2} \tilde G + \|\delta_i\|_*)^2}{2 \beta P_i (1+p_i) p_i } 
    + \frac{\langle \delta_i, u - u_{i-1} \rangle}{P_i (1 + p_i)} \right].
    \end{split}
\end{align}
$\tilde u_t$ is a convex combination of $\tilde u_{t-1}$ and  $u_t$ (Line \ref{tilde_u_t} of {\tt PS} procedure). In turn, $\tilde u_{t-1}$ is also a combination $\tilde u_{t-2}$ and  $u_{t-1}$. Continuing further, we have that $\tilde u_t$ is a convex combination of $u_t$, $u_{t-1}$, \ldots $u_1$. Using the definitions $\theta_t$ \eqref{theta_t_main} + \eqref{p_t_main} and $\tilde u_t$  (Line \ref{tilde_u_t} of {\tt PS} procedure) we have
\begin{align}
\label{temp34}
\tilde u_t &= (1-\theta_t) \Tilde u_{t-1} + \theta_t u_t
\notag\\
&= \left(1 - \frac{P_{t-1} - P_t}{(1 - P_t)P_{t-1}}\right) \tilde u_{t-1} + \frac{P_{t-1} - P_t}{(1 - P_t)P_{t-1}} u_t
\notag\\
&=\frac{P_t}{1 - P_t} \left(\frac{1 - P_{t-1}}{P_{t-1}}\tilde u_{t-1} + \frac{1}{P_t(1 + p_t)}u_t \right),\nonumber\\
&= \frac{P_t}{1 - P_t} \left(\frac{1 -P_{t-1}}{P_{t-1}}(1 - \theta_{t-1})\tilde u_{t-2} + \frac{1 - P_{t-1}}{P_{t-1}} \theta_{t-1} u_{t-1} + \frac{1}{P_t(1 + p_t)}u_t \right) \nonumber\\
&= \frac{P_t}{1 - P_t} \left(\frac{1 -P_{t-1}}{P_{t-1}}\Bigg(1 - \frac{P_{t-2} - P_{t-1}}{(1 - P_{t-1})P_{t-2}}\right)\tilde u_{t-2} 
\notag\\
&\hspace{2cm}+ \frac{1 - P_{t-1}}{P_{t-1}} \frac{P_{t-2} - P_{t-1}}{(1 - P_{t-1})P_{t-2}} u_{t-1} + \frac{1}{P_t(1 + p_t)}u_t \Bigg) \nonumber\\
&= \frac{P_t}{1 - P_t} \left(\frac{1 -P_{t-2}}{P_{t-2}}\tilde u_{t-2} + \frac{1}{P_{t-1} (1 + p_{t-1})} u_{t-1} + \frac{1}{P_t(1 + p_t)}u_t \right) \nonumber\\
&= \ldots = \sum\limits_{i=1}^t \frac{P_t}{P_i(1 + p_i)(1 - P_{t})} \cdot u_i.
\end{align}
\ab{Combining \eqref{temp34}, \eqref{temp21} and using convexity of $\Phi$, we get
\begin{align*}
    \frac{\beta}{1 - P_{t}} V(u_t, u) + &(\Phi(\tilde u_t)  - \Phi(u)) \notag\\
    \leq& \frac{\beta}{1 - P_{t}} V(u_t, u) + \sum_{i=1}^t \frac{P_{t}}{P_i (1+p_i) (1 - P_{t}) }\cdot (\Phi(u_i)  - \Phi(u)) 
    \notag\\
    \leq& \frac{\beta P_{t}}{1 - P_{t}} V(u_{0}, u) +  \frac{P_{t}}{1 - P_{t}}\sum_{i=1}^t \left[ \frac{(\abc{2} \tilde G + \|\delta_i\|_*)^2}{2 \beta P_i (1+p_i) p_i } 
    + \frac{\langle \delta_i, u - u_{i-1} \rangle}{P_i (1 + p_i)} \right]. 
\end{align*}
}
\EndProof

\begin{lemma}\label{first_theorem} Assume that $\{ p_t\}_{t\ge 1}$, $\{\theta_t\}_{t\ge 1}$, $\{\beta_k\}_{k \ge 1}$, $\{\gamma_k\}_{k\ge 1}$ in Algorithm 1 satisfy \eqref{theta_t_main}, \eqref{p_t_main} and
\begin{align}
    \label{gamma_main}
    &\gamma_1 = 1,~~~~\beta_k - L \gamma_k \geq 0,~~~ k\geq 1,
\\
&\Gamma_k = 
\begin{cases} 
            1, & k = 1,\\
            (1 - \gamma_k)\Gamma_{k-1}, & k > 1,
        \label{gamma_kk_main}
\end{cases}
\\
    \label{gamma_k_main_app}
    &\frac{\gamma_k \beta_k}{\Gamma_k(1 - P_{T_k})} \leq \frac{\gamma_{k-1} \beta_{k-1}}{\Gamma_{k-1}(1 - P_{T_{k-1}})} ,~~~ k\geq 2.
\end{align}
Then 
\begin{align} 
    \label{t2_1}
    \E[\Psi(\overline x_N) - \Psi(x^*)]
    \leq \frac{\Gamma_N \beta_1}{1 - P_{T_1}} V(x_0, x^*) 
    + \Gamma_N \sum\limits_{k=1}^N\sum\limits_{i=1}^{T_k} \frac{(\abc{2} \tilde G^2 + \rho^2)\gamma_k P_{T_k}}{\beta_k \Gamma_k(1 - P_{T_k})p_i^2 P_{i-1}},
\end{align}
where $P_t$ is from \eqref{p_t_main}, $\tilde G$ is from \eqref{tG}, $x^*$ is the solution for the problem \eqref{problem_orig}, 
\begin{align} \label{problem}
 \Psi(x) &\eqdef F(x) + g(x),
\\
\label{sig_main}
    \rho^2 &\eqdef 24 n p^2(n) G^2 + \frac{4n^2p^2(n)\sigma^2}{r^2}.
\end{align}
\end{lemma}
\textit{Proof of Lemma \ref{first_theorem}:}
The function $g$ is $L$-smooth. Hence, with \eqref{g-L-smooth},
we obtain:
\begin{align*}
g(\overline x_k) &\leq g(\underline{x}_k) + \langle \nabla g(\underline{x}_k), \overline x_k - \underline{x}_k \rangle + \frac{L}{2} \|\overline x_k - \underline{x}_k\|^2.
\end{align*}
Then we use $\overline x_k = (1 - \gamma_k)\overline x_{k-1} + \gamma_k \Tilde x_k$ and $\overline x_k - \underline{x}_k = \gamma_k (\tilde x_k - x_{k-1})$ (Lines \ref{under_bar_x_k} and \ref{bar_x_k} of Algorithm \ref{alg}) and get
\begin{align*}
    g(\overline x_k) \leq& g(\underline{x}_k) + \langle \nabla g(\underline{x}_k), \overline x_k - \underline{x}_k \rangle + \frac{L}{2} \|\overline x_k - \underline{x}_k\|^2 
    \nonumber\\
    =& g(\underline{x}_k) + \langle \nabla g(\underline{x}_k), (1 - \gamma_k)\overline x_{k-1} + \gamma_k \Tilde x_k - \underline{x}_k \rangle + \frac{L \gamma_k^2}{2} \|\tilde x_k - x_{k-1}\|^2
    \\
    =& g(\underline{x}_k) + (1 - \gamma_k) \langle \nabla g(\underline{x}_k), \overline x_{k-1} - \underline{x}_k \rangle + \gamma_k \langle \nabla g(\underline{x}_k), \Tilde x_k - \underline{x}_k \rangle 
    \\
    &+ \frac{L \gamma_k^2}{2} \|\tilde x_k - x_{k-1}\|^2
    \\
    =& (1-\gamma_k) \left( g(\underline{x}_k) + \langle \nabla g(\underline{x}_k), \overline x_{k-1} - \underline{x}_k \rangle \right) \\
    &+ \gamma_k (g(\underline{x}_k) + \langle \nabla g(\underline{x}_k), \tilde x_{k} - \underline{x}_k \rangle ) + \frac{L \gamma_k^2}{2} \|\tilde x_k - x_{k-1}\|^2\nonumber \\
    \leq& (1-\gamma_k) g(\overline x_{k-1}) + \gamma_k \left[ g(\underline{x}_k) + \langle \nabla g(\underline{x}_k), \tilde x_{k} - \underline{x}_k \rangle +\beta_k V(x_{k-1}, \tilde x_k) \right] \\
    &- \gamma_k \beta_k V(x_{k-1}, \tilde x_k) + \frac{L \gamma_k^2}{2} \|\tilde x_k - x_{k-1}\|^2\nonumber.
\end{align*}
Here we also used the convexity of $g$: $g(\underline{x}_k) + \langle \nabla g(\underline{x}_k), \overline x_{k-1} - \underline{x}_k \rangle \leq g(\overline x_{k-1})$.
By the strong convexity of $V$ (Lemma \ref{eq:bregman_key_property}) \ab{and \eqref{gamma_main}} we get: 
\begin{align*}
g(\overline x_k) \leq& (1-\gamma_k) g(\overline x_{k-1}) + \gamma_k \left[ g(\underline{x}_k) + \langle \nabla g(\underline{x}_k), \tilde x_{k} - \underline{x}_k \rangle +\beta_k V(x_{k-1}, \tilde x_k) \right] \nonumber\\
& - \left( \gamma_k \beta_k - L \gamma_k^2 \right) V(x_{k-1}, \tilde x_k) \nonumber \\
\leq& (1-\gamma_k) g(\overline x_{k-1}) + \gamma_k \left[ g(\underline{x}_k) + \langle \nabla g(\underline{x}_k), \tilde x_{k} - \underline{x}_k \rangle +\beta_k V(x_{k-1}, \tilde x_k) \right],
\end{align*}
Using the convexity of $F$ and Line \ref{bar_x_k} of Algorithm \ref{alg}, we obtain:
\begin{equation*}
F(\overline x_k) \leq (1-\gamma_k) F(\overline x_{k-1})  + \gamma_k F(\tilde x_k).
\end{equation*}
Summing up previous two inequalities, and using the definition \eqref{problem}, we have
\begin{align*}
\Psi(\overline x_k) 
\leq& (1-\gamma_k) \Psi(\overline x_{k-1}) + \gamma_k \left[ F(\tilde x_k) + g(\underline{x}_k) + \langle \nabla g(\underline{x}_k), \tilde x_{k} - \underline{x}_k \rangle +\beta_k V(x_{k-1}, \tilde x_k) \right]
\\
=&
(1-\gamma_k) \Psi(\overline x_{k-1}) + \gamma_k \left[ F(\tilde x_k) + h_k (\tilde x_{k}) +\beta_k V(x_{k-1}, \tilde x_k) \right]
\end{align*}
In the last step, we took into account the notation of $h_k$ from Line \ref{under_bar_x_k} of Algorithm \ref{alg}. 
\abc{Similarly to \eqref{Phi_main} we can introduce the definition:}
\ab{$\Phi_k(u) \eqdef h_k(u) + F(u) + \beta_k V(x_{k-1},u)$.} Then
\begin{equation}
\label{temp4_1}
\Psi(\overline x_k) - \Psi(u) \leq (1- \gamma_k) [\Psi(\overline x_{k-1})  - \Psi(u)] + \gamma_k [\Phi_k(\tilde x_k) - \Psi(u)].
\end{equation}
The convexity of the function $g$ gives
\begin{align*}
    \Phi_k(u) =& h_k(u) + F(u) + \beta_k V(x_{k-1},u)
    \notag\\
    =& g(x_k) + \langle \nabla g(x_k), u-x_k \rangle + F(u) + \beta_k V(x_{k-1},u)
    \notag\\
    \leq& g(u) + F(u) + \beta_k V(x_{k-1},u)
    \notag\\
    =&
    \Psi (u) + \beta_k V(x_{k-1},u).
\end{align*}
After small rearrangement we have
\begin{align}
\label{temp888}
    - \Psi (u) 
    \leq - \Phi_k(u) + \beta_k V(x_{k-1},u).
\end{align}
Combining \eqref{temp888} and \eqref{temp4_1}, we get
\begin{equation}
\label{temp4}
\Psi(\overline x_k) - \Psi(u) \leq (1- \gamma_k) [\Psi(\overline x_{k-1})  - \Psi(u)] + \gamma_k [\Phi_k(\tilde x_k) - \Phi_k(u)] + \gamma_k \beta_k V(x_{k-1},u).
\end{equation}
For the {\tt PS} subroutine with the input $(h_k, x_{k-1}, \beta_k, T_k)$ one can use the results of Lemma \ref{third_lemma} and get for all $u\in \X$
\begin{align}
        \label{1.65}
        \begin{split}
        \beta_k (1 - P_{T_k})^{-1} &V(x_k, u) + \left[\Phi_k(\tilde x_k) - \Phi_k(u)\right]\\
        \leq& \beta_k P_{T_k}(1 - P_{T_k})^{-1}V(x_{k-1},u)
        \\
        &+ \frac{P_{T_k}}{1 - P_{T_k}} \sum\limits_{i=1}^{T_k} \frac{1}{p_i P_{i-1}}  \left[\frac{(\abc{2} \tilde G + \| \delta_{k,i}\|_*)^2}{2 \beta_k p_i} + \langle \delta_{k,i}, u-u_{k,i-1} \rangle\right].
        \end{split}
\end{align}
Here we also used that $(x_k, \tilde x_k)$ is the output of {\tt PS}. 
Combing of \eqref{temp4} and \eqref{1.65} gives for all $u \in \X$:
\begin{align*}
\Psi(\overline x_k) - \Psi(u)  
\leq&  (1 - \gamma_k) [\Psi(\overline x_{k-1})  - \Psi(u)]
+ \gamma_k \Bigg\{
\frac{\beta_k}{1 - P_{T_k}} [V(x_{k-1}, u) - V(x_k, u)]  \nonumber\\
 & + \frac{P_{T_k}}{1 - P_{T_k} } 
\sum_{i=1}^{T_k} \frac{1}{p_i P_{i-1}} \left[\frac{\left(\abc{2} \tilde G + \|\delta_{k,i}\|_*\right)^2}{2\beta_k p_i } 
+ \langle \delta_{k,i}, u - u_{k,i-1} \rangle \right]\Bigg\}.
\end{align*}
Now we apply Lemma \ref{lem:lemma3.17_lan} with $w_k = \gamma_k$, $W_k = \Gamma_k$ and $\Delta_k = \Psi(\overline x_k) - \Psi(u) $ and get
\begin{align}
\Psi(\overline x_N) &- \Psi(u)   \notag\\
\leq &\Gamma_N (1- \gamma_1) [\Psi(\overline x_{0})  - \Psi(u)] \nonumber + \Gamma_N 
\sum_{k=1}^N \frac{\beta_k \gamma_k}{\Gamma_k (1 - P_{T_k})}
\left[
V(x_{k-1}, u) - V(x_k, u) \right] 
\nonumber \\
 & + \Gamma_N \sum_{k=1}^N 
 \frac{\gamma_k P_{T_k}}{\Gamma_k (1 - P_{T_k}) }\sum_{i=1}^{T_k} \frac{1}{p_i P_{i-1}} \left[\frac{\left(\abc{2} \tilde G + \|\delta_{k,i}\|_*\right)^2}{2 \beta_k p_i } 
+ \langle \delta_{k,i}, u - u_{k,i-1} \rangle 
\right]
\notag\\
= &\Gamma_N 
\sum_{k=1}^N \frac{\beta_k \gamma_k}{\Gamma_k (1 - P_{T_k})}
\left[
V(x_{k-1}, u) - V(x_k, u) \right] 
\nonumber \\
 & + \Gamma_N \sum_{k=1}^N 
 \frac{\gamma_k P_{T_k}}{\Gamma_k (1 - P_{T_k}) }\sum_{i=1}^{T_k} \frac{1}{p_i P_{i-1}} \left[\frac{\left(\abc{2} \tilde G + \|\delta_{k,i}\|_*\right)^2}{2 \beta_k p_i } 
+ \langle \delta_{k,i}, u - u_{k,i-1} \rangle 
\right]
.\label{eq:technical_sliding}
\end{align}
Here we also used that $\gamma_1 = 1$. From  \eqref{gamma_k_main_app} we obtain that
\begin{align*}
\sum_{k=1}^N \frac{\beta_k\gamma_k}{\Gamma_k (1 - P_{T_k})}
&\left[
V(x_{k-1}, u) -  V(x_k, u) \right] \notag\\
&\abc{\leq} \sum_{k=1}^N \left[\frac{\beta_{k}\gamma_{k}}{\Gamma_{k} (1 - P_{T_{k}})}
V(x_{k-1}, u) -  \frac{\beta_{k+1}\gamma_{k+1}}{\Gamma_{k+1} (1 - P_{T_{k+1}})} V(x_{k}, u)\right]
\\
&\abc{=} \frac{\beta_1 \gamma_1  }{\Gamma_1 (1 - P_{T_1})} V(x_0, u) - \frac{\beta_{N+1}  \gamma_{N+1} }{\Gamma_{N+1} (1 - P_{T_{N+1}})} V(x_{N+1}, u).
\end{align*}
With $\gamma_1 = \Gamma_1 = 1$, $P_{T_N} < 1$ and $V(x_N, u) \ge0$ we get
\begin{align}
\label{temp6}
\sum_{k=1}^N \frac{\beta_k\gamma_k}{\Gamma_k (1 - P_{T_k})} \left[
V(x_{k-1}, u) -  V(x_k, u) \right] &\leq \frac{\beta_1}{1 - P_{T_1}} V(x_0, u).
\end{align}
Substituting \eqref{temp6} into \eqref{eq:technical_sliding}, we get
    \begin{align*}
    \Psi(\overline x_N) - \Psi(u) \leq& \frac{\Gamma_N \beta_1}{1 - P_{T_1}} V(x_0,u) \notag\\
    &+
    \Gamma_N \sum\limits_{k=1}^N  \frac{\gamma_k P_{T_k}}{\Gamma_k(1 - P_{T_k})} \sum\limits_{i=1}^{T_k} \frac{1}{p_i P_{i-1}} \Bigg[ \frac{(\abc{2} \tilde G^2+ \|\delta_{k,i}\|_*^2)}{\beta_k p_i} + \langle \delta_{k,i}, u - u_{k,i-1}\rangle \Bigg].
    \end{align*}
\ab{
Then we substitute $u = x^*$ and take the full expectation
    \begin{align*}
    \EE[\Psi&(\overline x_N) - \Psi(x^*)]  \\
    \leq& \frac{\Gamma_N \beta_1 }{1 - P_{T_1}} V(x_0, x^*) \\ &+ 
    \Gamma_N \sum\limits_{k=1}^N  \frac{\gamma_k P_{T_k}}{\Gamma_k(1 - P_{T_k})} \sum\limits_{i=1}^{T_k} \frac{1}{p_i P_{i-1}} \Bigg[ \frac{(\abc{2} \tilde G^2+ \EE[\|\delta_{k,i}\|_*^2])}{\beta_k p_i} + \EE\left[\langle \delta_{k,i}, x^* - u_{k,i-1}\rangle\right] \Bigg].
    \end{align*}
Using the definition of $\delta_{k,i}$ from \eqref{delta_main}, one can obtain that $u_{k,i-1}$ does not depend on $\delta_{k,i}$. Therefore, by \eqref{apr_grad_diff}, we get
\begin{equation}
    \mathbb{E}[\langle\delta_{k,i}, x^* - u_{k,i-1}\rangle] = \mathbb{E}[\langle \E_{\delta_{k,i}}[\delta_{k,i}], x^* - u_{k,i-1}\rangle] = 0.
\end{equation}    
Whence we obtain
    \begin{align}
    \label{eq:24_88}
    \begin{split}
        \EE&\left[\Psi(\overline x_N) - \Psi(x^*)\right]  \\
    &\leq\frac{\Gamma_N \beta_1 }{1 - P_{T_1}} V(x_0, x^*) + 
    \Gamma_N \sum\limits_{k=1}^N  \frac{\gamma_k P_{T_k}}{\Gamma_k(1 - P_{T_k})} \sum\limits_{i=1}^{T_k} \frac{1}{p_i P_{i-1}} \Bigg[ \frac{(\abc{2} \tilde G^2+ \EE[\|\delta_{k,i}\|_*^2])}{\beta_k p_i} \Bigg].
    \end{split}
    \end{align}
Next, we estimate $\mathbb{E}[\|\delta_{k,i}\|^2_{*}]$: 
\begin{align*}
    \mathbb{E}[\|\delta_{k,i}\|^2_{*}] &= \mathbb{E}[\|\Tilde{f'_r}(u_{k,i-1}, \xi^{\pm}_{k,i-1}) - \nabla F(u_{k,i-1})\|^2_{*}] 
    \\
    &\leq  2\mathbb{E}\|\Tilde{f'_r}(u_{k,i-1}, \xi^{\pm}_{k,i-1}) \|^2_{*} + 2\mathbb{E}\|\nabla F(u_{k,i-1})\|^2_{*}.
\end{align*}
Using the results of Lemma \ref{first_lemma} and \ref{second_lemma}, we get
\begin{align*}
    \mathbb{E}[\|\delta_{k,i}\|^2_{*}] &\leq 2p^2 (n) \left(8n G^2 + \frac{2n^2\sigma^2}{r^2} \right) + 8 n p^2 (n) G^2. \notag
\end{align*}
It remains to substitute this estimate into \eqref{eq:24_88} and take into account \eqref{sig_main} to complete the proof.
}
    
\EndProof
Now we are ready to prove the main theorem.

\begin{theorem}[Theorem \ref{cor:main}]\label{cor:main_app} Suppose that $\{ p_t\}_{t\ge 1}$, $\{\theta_t\}_{t\ge 1}$ are
\begin{equation}
    \label{pt_main}
    p_t = \frac{t}{2}, ~~~ \theta_t = \frac{2(t+1)}{t(t+3)},
\end{equation}
$N$ is given, $\{\beta_k \}_{k \geq 1}$, $\{\gamma_k\}_{k \geq 1}$, $\{T_k\}_{k \geq 1}$ are
\begin{equation}
    \label{k_main}
    \beta_k = \frac{2L}{k},~~~\gamma_k = \frac{2}{k+1},~~~T_k = \max\left\{1; \frac{N(\abc{2} \tilde G^2 + \rho^2)k^2}{D_{\X,V}^2 L^2} \right\}
\end{equation}
for $\tilde G$ is  from \eqref{tG} and $\rho$ is from \eqref{sig_main}. Then for all $ N \geq 1$
\begin{equation}
    \label{col1_main}
    \mathbb{E}[\Psi_0(\overline x_N)- \Psi_0(x^*)] \leq\frac{20LD_{\X,V}^2}{N(N+1)} + \abc{2rG}.
\end{equation}
Additionally, the total number {\tt PS} procedure \abc{iterations} is
    \begin{equation}
        \label{original_total_main_app}
        T^{\text{total}} = \frac{(N+1)^4}{D_{\X,V}^2 L^2} \cdot \left( \abc{20} p^2(n) nG + \frac{2p^2(n) n^2\sigma^2}{r^2}\right) \abc{+ N}.
    \end{equation}
\end{theorem}
\textit{Proof of Theorem \ref{cor:main_app}:} First of all, we need to verify that our choice of parameters satisfies the conditions \eqref{theta_t_main}, \eqref{gamma_main}, \eqref{gamma_k_main_app}.
It is easy to see that for our $p_t$ 
\begin{equation}
    \label{Pt1}
    P_t = \frac{2}{(t+1)(t+2)}.
\end{equation}
One can check that such $P_t$ and $\theta_t$ from \eqref{pt_main} satisfy the condition \eqref{theta_t_main}.
Also with $T_k$ from \eqref{k_main}, we get
    \begin{equation}
    \label{Pt2}
    P_{T_k} \leq P_{T_{k-1}} \leq \ldots \leq P_{T_1} \leq \frac{1}{3}.
    \end{equation}
It is also easy to verify that for our $\gamma_k$ 
    \begin{equation}
    \label{gk}
    \Gamma_k = \frac{2}{k(k+1)}
    \end{equation}
Moreover, one can note that $\beta_k$ and $\gamma_k$ from \eqref{k_main} fit the inequality \eqref{gamma_main}. Finally, by \eqref{k_main}, \eqref{Pt1}, \eqref{Pt2}, \eqref{gk} we verify assumption \eqref{gamma_k_main_app}. 

Now, we are ready to prove \eqref{col1_main}. Simple calculations and relations \eqref{pt_main}, \eqref{Pt1} imply
    \begin{equation*}
    \sum\limits_{i=1}^{T_k} \frac{1}{p_i^2 P_{i-1}} = 2 \sum\limits_{i=1}^{T_k} \frac{i+1}{i} \leq 4 T_k.
    \end{equation*}
Next, from this estimate we can obtain
    \begin{equation*}
    \sum\limits_{i=1}^{T_k} \frac{\gamma_k P_{T_k}}{\Gamma_k \beta_k (1 - P_{T_k}) p_i^2 P_{i-1}} = \frac{\gamma_k P_{T_k}}{\Gamma_k \beta_k (1 - P_{T_k})} \sum\limits_{i=1}^{T_k} \frac{1}{p_i^2 P_{i-1}} \leq \frac{4 \gamma_k P_{T_k} T_k}{\Gamma_k \beta_k (1 - P_{T_k})}.
    \end{equation*}
Substituting $\gamma_k$, $\beta_k$ from \eqref{k_main} and $\Gamma_k$ from \eqref{gk}, one can obtain   
    \begin{equation*}
    \sum\limits_{i=1}^{T_k} \frac{\gamma_k P_{T_k}}{\Gamma_k \beta_k (1 - P_{T_k}) p_i^2 P_{i-1}} \leq \frac{4 \gamma_k P_{T_k} T_k}{\Gamma_k \beta_k (1 - P_{T_k})} \leq \frac{2 k^2 P_{T_k} T_k}{L (1 - P_{T_k})}.
    \end{equation*}
Using \eqref{Pt2}, we can note that $1 - P_{T_k} \geq \tfrac{2}{3}$. Also, substituting \eqref{Pt1} for $P_{T_k}$, we get
    \begin{equation}
    \label{p3}
    \sum\limits_{i=1}^{T_k} \frac{\gamma_k P_{T_k}}{\Gamma_k \beta_k (1 - P_{T_k}) p_i^2 P_{i-1}} \leq  \frac{2 k^2 P_{T_k} T_k}{L (1 - P_{T_k})} \leq \frac{3 k^2 P_{T_k} T_k}{L} =  \frac{6k^2 T_k}{L(T_k + 1)(T_k + 2)}.
    \end{equation}
Finally, we use the statement \eqref{t2_1} of Lemma \ref{first_theorem} and \eqref{p3}
\begin{align*} 
    \E[\Psi(\overline x_N) - \Psi(x^*)]
    \leq& \frac{\Gamma_N \beta_1}{1 - P_{T_1}} V(x_0, x^*) 
    + 6\Gamma_N (\abc{2} \tilde G^2 + \rho^2) \sum\limits_{k=1}^N \frac{k^2 T_k}{L(T_k + 1)(T_k + 2)}
    \\
    \leq& \frac{\Gamma_N \beta_1}{1 - P_{T_1}} V(x_0, x^*) 
    + 6\Gamma_N (\abc{2} \tilde G^2 + \rho^2) \sum\limits_{k=1}^N \frac{k^2}{L T_k}.
\end{align*}
Substituting $T_k$ from \eqref{k_main}, we have
\begin{align*} 
    \E[\Psi(\overline x_N) - \Psi(x^*)]
    &\leq \frac{\Gamma_N \beta_1}{1 - P_{T_1}} V(x_0, x^*) 
    + 6\Gamma_N \sum\limits_{k=1}^N \frac{D_{\X,V}^2 L}{ N} \\
    &\leq\Gamma_N \left( \frac{ \beta_1}{1 - P_{T_1}} V(x_0, x^*) 
    + 6 D_{\X,V}^2 L\right).
\end{align*}
It remains to put $\Gamma_N$ from \eqref{gk}, $\beta_1$ from \eqref{k_main} and $1 - P_{T_1} \geq \tfrac{2}{3}$. 
    \begin{align*}
    \E[\Psi(\overline x_N) - \Psi(x^*)] \leq \frac{2L}{N(N+1)}(3V(x_0, x^{*}) + 6D_{\X,V}^2) \leq \frac{20LD_{\X,V}^2}{N(N+1)}.
    \end{align*} 
To prove \eqref{col1_main} it remains to estimate the relationship between $\Psi_0$ and $\Psi$. In particular, one can note that for all $x \in \X$ we have that $|\Psi_0(x) - \Psi(x)| = |f(x) - F(x)|$. This term can be bounded by \eqref{orig_sm_main_app}. Therefore, we obtain
\begin{align*}
    \E[\Psi_0(\overline x_N) - \Psi_0(x^*)] 
    &\leq \E[|\Psi_0(\overline x_N) - \Psi(\overline x_N)| + |\Psi_0(x^*) - \Psi(x^*)| + \Psi(\overline x_N) - \Psi(x^*)] 
    \\
    &\leq 2rG + \frac{20LD_{\X,V}^2}{N(N+1)}.
\end{align*} 
To prove the bounds \eqref{original_total_main_app} we use \eqref{k_main} for $T_k$
 \begin{align*}
    T^{\text{total}} &=
    \sum\limits_{i=1}^{N} T_k 
    \\
    &{\leq}  \sum\limits_{i=1}^{N} \left(\frac{N(\abc{2} \tilde G^2 + \rho^2)k^2}{D_{\X,V}^2L^2} + 1\right) 
    \nonumber\\
    &= \frac{1}{6} \frac{N^2(N+1)(2N+1)(\abc{2} \tilde G^2 + \rho^2)}{D_{\X,V}^2 L^2} + N \notag\\
    &\leq \frac{1}{3} \frac{(N+1)^4(\abc{2} \tilde G^2 + \rho^2)}{D_{\X,V}^2 L^2} + N. 
\end{align*} 
Substituting $\tilde G$ and $\rho$ from \eqref{tG} and \eqref{sig_main} gives the final result.
\EndProof

\abc{
\section{Proof of Lemma \ref{lem:prop}}

\begin{lemma}
Let the functions $f_m$ from \eqref{centr} and \eqref{decentr}  be $G$-Lipschitz w.r.t. $\| \cdot \|_2$, the local noises $\xi_m$ for all $f_m(x_m, \xi_m) = f_m(x_m) + \xi_m$ be independent, unbiased and bounded: $\E \xi_m = 0$, $\E[\xi^2_m] \leq \sigma^2$, the Euclidean diameter of the set $\X$ be equal to $D_{\X}$. Then the following facts for \eqref{centr} and \eqref{decentr} are valid:
\begin{itemize}
    \item $f(x_1,\ldots,x_M)$ is $(G/\sqrt{M})$-Lipschitz w.r.t. $\| \cdot \|_2$ of $\textbf{x} = [x_1^T, \ldots x_M^T]^T$;
    \item $g(x_1, \ldots, x_M)$ from \eqref{centr} is $(\lambda /M)$-smooth w.r.t. $\| \cdot \|_2$ of $\textbf{x}$, and $g$ from \eqref{decentr} is $(\lambda \lambda_{\max}(W)/M)$-smooth w.r.t. $\| \cdot \|_2$ of $\textbf{x}$;
    \item  the Euclidean diameter $D_{\X^M}$ of the set $\X^M$ is equal to $\sqrt{M} D_{\X}$;
    \item the noise of the function $f(x_1,\ldots,x_M, \xi_1, \ldots, \xi_M)$ is unbiased and bounded by $\sigma^2/M$.
\end{itemize}
\end{lemma}
\textit{Proof of Lemma~\ref{lem:prop}:}
We start from the first point:
\begin{align*}
   |f(x_1,\ldots,x_M) - f(y_1,\ldots,y_M) | &= \left| \frac{1}{M}\sum\limits_{m=1}^M [f_m(x_m) - f_m(y_m)]\right| 
   \\
   &\leq \frac{1}{M}\sum_{m=1}^M |f_m(x_m) - f_m(y_m)| 
   \\
   &\leq \frac{G}{M}\sum_{m=1}^M \|x_m - y_m\|_2. 
\end{align*}
Here we used the $G$-Lipschitzness of $f_m$. Next, we square this expression and apply \eqref{eq:squared_sum} to get
\begin{align*}
|f(x_1,\ldots,x_M) - f(y_1,\ldots,y_M)|^2 
&\leq \frac{G^2}{M^2} \left(\sum_{m=1}^M \|x_m - y_m\|_2\right)^2 
\\
&\leq \frac{G^2}{M} \sum_{m=1}^M \|x_m - y_m\|^2_2 
\\
&= \frac{G^2}{M} \| \textbf{x} - \textbf{y} \|_2^2.
\end{align*}
It finishes the proof that $f(x_1,\ldots,x_M)$ is $G/\sqrt{M}$--Lipschitz. 

We move to the smoothness of $g$ and consider \eqref{centr}. One can find $\nabla g$ from \eqref{centr} as follows
\begin{align*}
    \nabla_{x_j}\bigg[\frac{\lambda}{2 M}\sum\limits_{m=1}^M &\|x_m - \bar x\|^2_2 \bigg]
    \\
    &= \nabla_{x_j}\left[\frac{\lambda}{2 M} \left\|x_j - \frac{1}{M} x_j - c \right\|^2_2 \right] + \nabla_{x_j}\left[\frac{\lambda}{2 M}\sum\limits_{m \neq j} \left\| \frac{1}{M} x_j  - \frac{d_m}{M}\right\|^2_2 \right] 
    \\
    &= \frac{\lambda}{M} \cdot \frac{M-1}{M} \left( \frac{M-1}{M} x_j - c \right) + \frac{\lambda}{M} \cdot \frac{1}{M^2} \sum\limits_{m \neq j} (x_j - d_m),
\end{align*}
where $c = \frac{1}{M}\sum\limits_{m \neq j} x_m$, $d_m = M \cdot x_m - \sum_{l \neq j} x_l$. By Definition \ref{def:smooth}, we obtain
\begin{align*}
\Bigg\| \nabla_{x_j}\left[\frac{\lambda}{2M}\sum\limits_{m=1}^M \|x_m - \bar x\|^2_2 \right] - \nabla_{y_j} & \left[\frac{\lambda}{2M}\sum\limits_{m=1}^M \|y_m - \bar y\|^2_2 \right] \Bigg\|^2_2 
\\
&= \left\| \lambda \cdot \frac{(M-1)^2}{M^3} (x_j - y_j) + \lambda \cdot \frac{M-1}{M^3} (x_j - y_j)\right\|^2_2
\\
&\leq \frac{\lambda^2}{M^2} \| x_j - y_j \|_2^2. 
\end{align*}
Summing over all $j$ from $1$ to $M$, we get
\begin{align*}
\Bigg\| \nabla_{x} & \left[\frac{\lambda}{2M}\sum\limits_{m=1}^M \|x_m - \bar x\|^2_2 \right] - \nabla_{y}  \left[\frac{\lambda}{2M}\sum\limits_{m=1}^M \|y_m - \bar y\|^2_2 \right] \Bigg\|^2_2 
\\
&= \sum_{j=1}^M \Bigg\| \nabla_{x_j}\left[\frac{\lambda}{2M}\sum\limits_{m=1}^M \|x_m - \bar x\|^2_2 \right] - \nabla_{y_j} \left[\frac{\lambda}{2M}\sum\limits_{m=1}^M \|y_m - \bar y\|^2_2 \right] \Bigg\|^2_2
\\
&\leq \frac{\lambda^2}{M^2} \sum_{j=1}^M \| x_j - y_j \|_2^2
\\
&= \frac{\lambda^2}{M^2} \| \textbf{x} - \textbf{y} \|_2^2. 
\end{align*}
It means that $g$ from \eqref{centr} is $\lambda/M$-smooth. Next, we consider $g$ from \eqref{decentr}. One can note that $\| \sqrt{W} X\|_F^2$ can be rewritten as follows: $\| (\sqrt{W} \otimes I) \textbf{x}\|_2^2$, where $I$ is the identity matrix of size $n$. It is easy to see that 
\begin{align*}
    \nabla_{\textbf{x}} \left[ \frac{\lambda}{2M}\| (\sqrt{W} \otimes I) \textbf{x}\|_2^2 \right] 
    &= \frac{\lambda}{M} \left[ (\sqrt{W} \otimes I)^T (\sqrt{W} \otimes I) \right] \textbf{x} 
    \\
    &= \frac{\lambda}{M} \left[ (\sqrt{W}^T \otimes I) (\sqrt{W} \otimes I) \right] \textbf{x} 
    \\
    &= \frac{\lambda}{M} \left[ \sqrt{W}^T \sqrt{W} \otimes I \right] \textbf{x}.
\end{align*}
Then we can estimate the smoothness constant:
\begin{align*}
\left\| \frac{\lambda}{M} \left[ \sqrt{W}^T \sqrt{W} \otimes I \right] \textbf{x} - \frac{\lambda}{M} \left[ \sqrt{W}^T \sqrt{W} \otimes I \right] \textbf{y} \right\|_2
&\leq 
\frac{\lambda}{M} \left\| \sqrt{W}^T \sqrt{W} \otimes I \right\|_2 \cdot \| \textbf{x} - \textbf{y}\|_2 
\\
&\leq 
\frac{\lambda}{M} \cdot \lambda_{\max}(\sqrt{W}^T \sqrt{W} \otimes I) \| \textbf{x} - \textbf{y}\|_2  
\\
&\leq 
\frac{\lambda}{M} \cdot \lambda_{\max}(\sqrt{W}^T \sqrt{W}) \| \textbf{x} - \textbf{y}\|_2 
\\
&\leq 
\frac{\lambda \lambda_{\max} (W)}{M} \| \textbf{x} - \textbf{y} \|_2. 
\end{align*}
Here we took into account the matrix analysis.

We work with the third bullet:
$$
\max_{(x_1, \ldots, x_M), (y_1, \ldots, y_M) \in \mathcal{X}^M }\| \textbf{x} - \textbf{y}\|^2_2 = \sum\limits_{m=1}^M \max_{x_m, y_m \in \X}\| x_m - y_m \|^2_2 \leq M D_{\X}^2. 
$$

Finally, we consider the point on the stochastic noise:
$$
\mathbb{E}\left[
\left(\frac{1}{M}\sum\limits_{m=1}^M \xi_m\right)^2 \right] =  \frac{1}{M^2}\sum\limits_{m=1}^M \mathbb{E}\left[\xi_m^2 \right] + \frac{1}{M^2}\sum\limits_{m \neq m'} \mathbb{E}\left[\langle \xi_m, \xi_{m'}\rangle \right] \leq \frac{\sigma^2}{M}. 
$$
Here we took into account that the independence of $\xi_m$ and $\xi_{m'}$, as well as their unbiasedness and boundedness to put $\mathbb{E}\left[\langle \xi_m, \xi_{m'}\right] = 0$ and to estimate $\mathbb{E}\left[\xi_m^2 \right]$ by $\sigma^2$. 

It finished the whole proof. 
\EndProof
}
\end{document}